\documentclass[11pt,twoside]{article}
\usepackage{amsmath,amsthm,amssymb,amsbsy,fancyhdr,graphicx,
            calc,ifthen,float,psfrag,mathrsfs,hhline}
\input epsf
\PassOptionsToPackage{ctagsplt,Righttag}{amstex}
\usepackage[colorlinks=true,citecolor=blue,linkcolor=blue,bookmarksnumbered=true,pdfstartview={FitH},hyperindex=true,linktoc=all,linktocpage=true]{hyperref}
\usepackage{mathtools}
\usepackage{framed}
\usepackage{dsfont}
\usepackage{ulem}
\usepackage{color}
\definecolor{darkgreen}{rgb}{0,0.65,0}
\definecolor{darkred}{rgb}{0.65,0,0}
\newcommand{\by}{\mathbf{y}}

\newcommand{\bnabla}{\boldsymbol{\nabla}}

\newcommand{\tr}{\operatorname{tr}}
\newcommand{\assembly}{\underset{e}{\textbf{\Large\textsf{A}\normalsize}}}
\newcommand{\secref}[1]{Section~\ref{#1}}
\newcounter{hours}
\newcounter{minutes}
\newcommand{\printtime}{\setcounter{hours}{\time/60}%
                        \setcounter{minutes}{\time-\value{hours}*60}%
\ifthenelse{\value{hours}<10}{0}{}\thehours:%
\ifthenelse{\value{minutes}<10}{0}{}\theminutes}
\def\lsp{\def\baselinestretch{0.75}}
\def\ssp{\def\baselinestretch{1.0}}
\def\dsp{\def\baselinestretch{1.37}}
\advance\textheight by \topskip
\usepackage{amsmath}
\usepackage{amssymb}
\usepackage{bm}
\usepackage{mathtools}
\usepackage[font=footnotesize]{caption}
\usepackage[font=footnotesize]{subcaption}
\usepackage{tabstackengine}
\usepackage{float}
\usepackage[ruled,vlined,linesnumbered,noend]{algorithm2e}
\SetAlgoNoLine 

\newcommand{\norm}[1]{\ensuremath{\left\| #1 \right\|}}

\newcommand{\pder}[2]{\ensuremath{\frac{\partial #1}{\partial #2}}} 
\newcommand{\pderTwo}[3]{\ensuremath{\frac{\partial^2 #1}{\partial #2 \partial #3}}} 
\newcommand{\pderH}[3]{\ensuremath{\frac{\partial^{#3} #1}{\partial #2^{#3}}}}

\newcommand{\Cbb}{\ensuremath{\mathbb{C}}}

\newcommand{\Rbb}{\ensuremath{\mathbb{R} }}

\newcommand\Abm{{\ensuremath{\bm{A}}}}
\newcommand\Bbm{{\ensuremath{\bm{B}}}}

\newcommand\Fbm{{\ensuremath{\bm{F}}}}

\newcommand\Ibm{{\ensuremath{\bm{I}}}}

\newcommand\Kbm{{\ensuremath{\bm{K}}}}

\newcommand\Mbm{{\ensuremath{\bm{M}}}}
\newcommand\Nbm{{\ensuremath{\bm{N}}}}

\newcommand\Qbm{{\ensuremath{\bm{Q}}}}
\newcommand\Rbm{{\ensuremath{\bm{R}}}}

\newcommand\Tbm{{\ensuremath{\bm{T}}}}

\newcommand\Xbm{{\ensuremath{\bm{X}}}}

\newcommand\abm{{\ensuremath{\bm{a}}}}
\newcommand\bbm{{\ensuremath{\bm{b}}}}

\newcommand\dbm{{\ensuremath{\bm{d}}}}

\newcommand\pbm{{\ensuremath{\bm{p}}}}
\newcommand\qbm{{\ensuremath{\bm{q}}}}
\newcommand\rbm{{\ensuremath{\bm{r}}}}

\newcommand\ubm{{\ensuremath{\bm{u}}}}
\newcommand\vbm{{\ensuremath{\bm{v}}}}

\newcommand\ybm{{\ensuremath{\bm{y}}}}
\newcommand\zbm{{\ensuremath{\bm{z}}}}

\newcommand\ubold{\ensuremath{\mathbf{u}}}

\newcommand\epsbold{{\ensuremath{\boldsymbol{\epsilon}}}}

\newcommand\taubold{{\ensuremath{\boldsymbol{\tau}}}}

\newcommand\thetabold{{\ensuremath{\boldsymbol{\theta}}}}
\newcommand\sigmabold{{\ensuremath{\boldsymbol{\sigma}}}}

\newcommand\xibold{{\ensuremath{\boldsymbol{\xi}}}}

\newcommand\zerobold{\ensuremath{\mathbf{0}}}

\newcommand{\dt}[1]{\Delta t_{#1}}
\newcommand\compactcdot{\mathinner{\mkern-2mu\cdotp\mkern-2mu}}

\begin{document}
\flushbottom
\pagestyle{empty}
\pagenumbering{arabic}
\ssp 
\title{\vspace{-0.75in}\bf 
  Optimization of process parameters in additive manufacturing based on the finite element method} 
\vspace{0.1in}
\author{\hspace{-0.25in} Jingyi WANG and Panayiotis 
PAPADOPOULOS\footnote{panos@berkeley.edu} 
\\[0.05in]
\sl\small Department of Mechanical Engineering, University of California, 
Berkeley, CA, USA 
\rm \normalsize}
\date{}
\maketitle
\small
\lsp
\par\noindent
\protect\vspace{0.1in}
\renewcommand{\contentsname}{\normalsize\centerline{Table of contents}}
\tableofcontents
\ssp
\protect\vspace{0.3in}
\normalsize
\begin{abstract}
A design optimization framework for process parameters of additive 
manufacturing based on finite element simulation is proposed. The finite element method 
uses a coupled thermomechanical model developed for fused deposition modeling from the authors' previous work. 
Both gradient-based and gradient-free optimization methods 
are proposed. The gradient-based approach, which solves a PDE-constrained optimization problem, 
requires sensitivities computed from the fully discretized finite element model. 
We show the derivation of the sensitivities and apply them in a projected gradient descent algorithm. 
For the gradient-free approach, we propose two distinct algorithms:
a local search algorithm called the method of local variations and a Bayesian optimization algorithm using Gaussian processes.  
To illustrate the effectiveness and differences of the methods, we provide 
two-dimensional design optimization examples using all three proposed algorithms. 
\end{abstract}
\par\noindent
\protect\vspace{0.1in}
\nopagebreak[5]
\begin{description}
\item[Keywords:] Additive manufacturing; sensitivity;
optimization; Bayesian optimization; 
\end{description}
\dsp
\newpage
\pagestyle{fancyplain}
\normalsize
%
%
\section{Introduction}\label{se:intro}
\par\noindent
Additive manufacturing (AM) has enjoyed substantial success in creating parts with complex geometries, 
shortening the product design time and lowering costs~\cite{mohamed2015,gibson2014,chua2010,upcraft2003,mansour2003,hopkinson2006}. 
Many additive manufacturing technologies have been developed in recent years  
including fused deposition modeling (FDM), selective laser sintering, {\em 
etc}. AM is adopted extensively in prototyping and is increasingly deployed to industrial production of parts for the medical, aerospace and 
automotive industries~\cite{thavornyutikarn2014,ahn2002}. 
The computer simulation of the AM processes is a topic of huge interest throughout the industry and academia~\cite{francois2017modeling,llnl-am-overview,overview2017,hajializadeh2019,kalita2003,Schoinochoritis2017,Jayanath2018}, where finite element method is one of the primary numerical tools used.
\par
FDM simulation models typically solves the transient heat transfer problem during the deposition, followed by the  
mechanics problem based on the temperature history of a material 
with temperature-dependent constitutive response~\cite{llnl-am-overview,overview2017,Michaleris2014,YZhang2006,patil2013,ding2014,Montevecchi2016}.
The process of deposition can be emulated using an ``active/inactive''  
element addition approach~\cite{Michaleris2014,michaleris2,michaleris3,Roberts2009,Yang2016,Montevecchi2016}. 
Heat transfer through conduction, convection and radiation are generally 
considered in simulation and the dependence of the thermal conductivity and 
heat capacity on temperature can be explicitly accounted for~\cite{Kolossov2004}. 
Previously, the authors proposed a fully coupled thermomechanical model for FDM~\cite{wang2021}, where the displacement and temperature fields are solved simultaneously and the deposition of new material is enabled through 
the creation of new elements. By accounting for existing displacement at each time step, 
the model can predict the quantities of interest such as maximum deformation, shape error, {\em etc}. 
\par
While AM simulation has progressed steadily, challenges remain in developing optimization tools for it. 
The optimization of process 
parameters of AM to achieve various design goals such as better surface finish quality, lower residual stress 
and less shape error for assembly is of great interest 
and importance~\cite{masood1996,groza2007,bellehumeur2004,chacon2017,YZhang2006,ding2014}. 
The process parameters to consider include chamber temperature, part orientation, printing speed, filament diameter, nozzle size and layer thickness~\cite{nancharaiah2010,horvath2007,wang2007,thrimurthulu2004,mukherjee2017}. 
The existing optimization tools for additive manufacturing are mostly simple, 
data-driven models using experimental data~\cite{mohamed2015,foroozmehr2016,west2001,Strano2013,Onwubolu2014,espin2015}.
Commonly used experiments include tensile, compression tests for the printed samples~\cite{ahn2002,ang2006,sood2010,percoco2012,rayegani2014,masood2010}, while dynamic mechanical tests are also performed in~\cite{arivazhagan2011,arivazhagan2012}.
Some of the techniques used for optimization are the Taguchi method~\cite{anitha2001,sood2009,nancharaiah2010,lee2005}, particle swarm optimization~\cite{sood2012} and artificial neural network (NN)~\cite{cho2000}. 
\par
In addition to experiments, an accurate finite element simulation is another natural candidate for optimization due  
to its flexibility and efficiency, 
and has been adopted by some~\cite{syrlybayev2021}. 
Process mapping of the parameters can be developed using finite element simulation and subsequently applied to optimization~\cite{Gockel2014,Krol2012,vasinonta2000,yin2022}. The finite element model used is often an uncoupled thermomechanical one with a fixed Lagrangian mesh~\cite{vastola2016,nickel2001}, which could raise the difficulty of sensitivity calculation for gradient-based optimization algorithms.
The optimization variables  of these studies include print speed, extrusion temperature and layer thickness~\cite{alafaghani2017,zhang2006}.
A purely geometric model is used in~\cite{morgan2016}.
\par
Bayesian optimization has been applied to many engineering design problems such as 
inverse problems~\cite{wang2004}, structural design~\cite{mathern2021} and 
robotics~\cite{calandra2016} with advanced technologies such as multi-fidelity 
surrogate models and independent constraints~\cite{bernardo2011,zuluaga2013}. 
Thus far, it has  been applied to a limited extent in the 
optimization of AM problems~\cite{zhang2021}, particularly for process parameters. Most of the approaches taken are 
experiment-driven or rely on geometric models, \textit{e.g.}, part orientation 
optimization of AM in~\cite{goguelin2021}. 
In~\cite{shapre2018,xue2020,hertlein2020},  machine learning and Bayesian optimization are used for the optimization of lattice structure for metamaterials, where the stiffness of the metamaterial is computed through finite element software.
A Bayesian optimization approach is adopted for the metal AM melt pool geometry optimization in~\cite{modal2020}. 
A data-driven Bayesian optimization method using finite element software to 
generate sample points is proposed in~\cite{xiong2019}. 
While~\cite{baturynska2018} argued 
for a conceptual framework to integrate experiments, finite element and machine 
learning together for the optimization of AM processes, the authors believe it remains a work in progress.
\par
In most of the existing optimization methods discussed above, the simulation and optimization are largely separated.
Additionally, process parameters as well as defects and uncertainty caused by the printing itself are often ignored. 
Further, the optimization algorithm applied might not have been systematically established and thus relies on exhaustive search. 
In this paper, we propose two approaches to contribute to the 
optimization workflow of AM process parameters, each with its own advantages. 
The first approach is a gradient-based optimization method, where we optimize the design objective 
constrained by partial differential equations (PDEs). 
A coupled finite element simulation model for FDM developed by the 
authors~\cite{wang2021} is used, whose fully discretized form is parameterized by the optimization variables, \textit{i.e.}, the process parameters. Then, we apply gradient-based optimization techniques to solve the optimization problem using the fully discretized sensitivities computed from the finite element model. 
The second approach is gradient-free and more suitable for problems where sensitivities are not easily available. 
We propose two algorithms that are differentiated by whether a surrogate model is used to approximate the objective, at least locally. 
The first algorithm, the method of local variations, advances through local function evaluations and step size updates. 
The second algorithm is a Bayesian optimization one that uses Gaussian process as its surrogate model, which is updated with new simulation data. 

\par
The organization of this paper is as follows. \secref{se:goveqn} contains a 
summary of the coupled thermomechanical finite element model. 
The optimization formulation is 
given in \secref{se:opt}. In~\secref{se:grad}, the gradient-based optimization 
algorithm is described, with derivations of sensitivities. 
In~\secref{se:grad-free}, we propose our gradient-free optimization algorithms
and discuss algorithmic parameter choices.  
Numerical examples of the optimization methods are presented in 
\secref{se:exp}. Finally, conclusions are included in \secref{se:conclusion}. 
%
%
\section{Review of continuum theory and finite element modeling}\label{se:goveqn}
\subsection{Continuum theory}\label{se:mechanics}
\par
In this paper, it is assumed that the printed body can be adequately modeled 
as an isotropic and homogeneous thermoelastic continuum undergoing 
infinitesimal deformation relative to its evolving reference configuration 
in the presence of large changes in its temperature. The material
comprising the body is locally endowed with a Helmholtz free energy 
function~${\psi=\bar\psi(\epsbold,\theta)}$ per unit volume, where~$\epsbold$ 
and~$\theta$ are the infinitesimal strain tensor and temperature, respectively. 
\par
The local form of linear momentum balance may be expressed as 
\begin{equation}\label{eqn:LMB_AMB}
   \rho_0\abm\ =\ \bnabla \compactcdot\sigmabold + \rho_0\bbm\ , 
\end{equation}
where~$\rho_0$ is the mass density per unit referential volume, 
$\abm$ is the acceleration, $\sigmabold$ is the stress tensor of the
infinitesimal theory, and~$\bbm$ is the body force per unit mass, 
while~``$\bnabla\compactcdot$'' denotes the divergence operator relative to the 
referential coordinates. The balance of angular momentum ensures that the 
stress tensor~$\sigmabold$ is symmetric. Moreover, 
the balance of energy takes the form 
\begin{equation}\label{eqn:energy}
\rho_0\dot e\ =\ 
\rho_0 r-\bnabla\compactcdot\qbm_0+\sigmabold\compactcdot\dot\epsbold\ ,
\end{equation}
in terms of the internal energy per unit mass~$e$, the heat supply 
per unit mass~$r$, and the referential heat flux 
vector~$\qbm_0=\bar \qbm_0 \left(\epsbold,\theta,\nabla\theta\right)$. 
\par
Upon invoking the Clausius-Duhem inequality, a standard procedure 
leads to the relations 
\begin{equation}\label{eqn:sh}
\sigmabold\ =\ \pder{\bar\psi\left(\epsbold,\theta
\right)}{\epsbold}\quad,\quad
\qbm_0\cdot\nabla\theta \ \leqslant\ 0\ , 
\end{equation}
as well as to the reformulation of the energy equation~\eqref{eqn:energy} as 
\begin{equation}\label{eqn:energy_bal}  
  c \ \dot\theta\ =\ -\bnabla\compactcdot\qbm_0 + \rho_0 r + 
                  \theta\Mbm \compactcdot\dot\epsbold\ ,
\end{equation}
where~$c=-\theta \pderH{\bar{\psi} \left(\epsbold,\theta\right)}{\theta}{2}$ is 
the heat capacity and
$\Mbm=\pderTwo{\bar{\psi}\left(\epsbold,\theta\right)}{\epsbold}{\theta}$
is the stress-temperature modulus, see, {\em e.g.}, \cite{wang2021} 
for a full derivation. 
\par
A Helmholtz free-energy function that 
can adequately represent materials under infinitesimal deformation and  
a finite temperature range~\cite{govindjee} is chosen according to 
\begin{equation}\label{const-helm-energy}
  \bar\psi(\epsbold,\theta)\ =\ 
\frac{1}{2} \epsbold\compactcdot\Cbb(\theta)\epsbold  
- \kappa(\theta)\ln(1+\tr\epsbold)\alpha(\theta-\theta_0)
+ \bar c
\left(\theta-\theta_0-\theta\ln{\frac{\theta}{\theta_0}}\right)\ ,
\end{equation}
where~$\alpha$ is the (constant) coefficient of thermal expansion,
$\theta_0$ is the reference temperature, and~$\bar c$ is the specific heat 
at the reference temperature. Also, the temperature-dependent isotropic elastic 
modulus in~\eqref{const-helm-energy} is expressed 
as~$\Cbb(\theta)=\tilde{f}\left(\theta \right)\Cbb_{0,m}$, 
where~$\tilde{f}(\cdot)$ is a smooth and positive function of temperature and
$\Cbb_{0,m}$ is the elastic modulus at a given temperature~$\theta_{0,m}$, not
necessarily equal to the reference temperature $\theta_0$. Likewise, 
the temperature-dependent bulk modulus in~\eqref{const-helm-energy} is  
approximated as~${\kappa(\theta)=\tilde{f}\left(\theta \right)\kappa_{0,m}}$, 
with~$\kappa_{0,m}$ being the bulk modulus at~$\theta_{0,m}$. 
The temperature function is defined here as 
\begin{equation}\label{const-helm-energy-abs-simp}
 \tilde{f}(\theta)\ =\ b\left(\frac{\theta}{\theta_{0,m}}\right)^a 
       + b\left(a-1\right) +
\left(1-ab\right)\frac{\theta}{\theta_{0,m}}\ , 
\end{equation}
where~$a$ and~$b$ are material constants and~$\tilde{f}(\theta_{0,m})=1$.
\par
Consistent with the restriction in~\eqref{eqn:sh}$_2$, the heat flux 
follows the isotropic form of Fourier's law, according to which 
\begin{equation}\label{eqn:fourier_law}
  \qbm_0\ =\ - k\nabla\theta\ ,
\end{equation}
where~$k$ is the (constant) thermal conductivity coefficient. Also, 
heat convection is included in the model, such that the flux of heat through 
the exterior boundary is given by 
\begin{equation}\label{eqn:convection-2}
 \centering
  \bar q\ =\ h\left(\theta_{\infty}-\theta \right)\ ,
\end{equation}
where~$h$ is the (constant) convection coefficient and~$\theta_{\infty}$ is 
the ambient temperature. 
\par
In view of the constitutive equation in~\eqref{const-helm-energy} and
upon taking into account~\eqref{eqn:sh}$_1$ and the definitions 
of the heat capacity $c$ and the stress-temperature modulus $\Mbm$, it follows that 
\begin{equation}\label{const-elastic-stress}
 \centering
 \begin{aligned}
 \sigmabold\ &=\ \tilde{f}\left(\theta \right)\left[\Cbb_{0,m}\epsbold - 
 \kappa_{0,m}\alpha\frac{\theta-\theta_0}{1+\tr\epsbold}\Ibm\right]\ , \\
 \Mbm\ &=\ \pder{\tilde{f}(\theta)}{\theta}\Cbb_{0,m}\epsbold
-\left[\pder{\tilde{f}(\theta)}{\theta} (\theta-\theta_0)+ 
\tilde{f}(\theta)\right] \kappa_{0,m}\alpha\frac{1}{1+\tr\epsbold}\Ibm\ , \\
  c\ &=\ -\theta\pderH{\tilde{f}(\theta)}{\theta}{2}\frac{1}{2}
 \epsbold\compactcdot\Cbb_{0,m}\epsbold +
\left[\theta\pderH{\tilde{f}(\theta)}{\theta}{2}(\theta-\theta_0)+\theta\pder{\tilde{f}(\theta)}{\theta}\right]\kappa_{0,m}\alpha \ln(1+\tr\epsbold) + \bar c\ ,
 \end{aligned}
\end{equation}
where~$\Ibm$ is the second-order identity tensor. 
%
\subsection{Finite element modeling}\label{se:fem}
\par
Let the printed body occupy the region~$\Omega_0$ in a reference configuration. 
The referential boundary~$\partial\Omega_0$ of the body is decomposed 
into Dirichlet parts~$\Gamma_{D,0}^u$, $\Gamma_{D,0}^t$ and 
Neumann parts~$\Gamma_{N,0}^u$, $\Gamma_{N,0}^t$ for displacement 
and temperature, respectively, such that 
$\overline{\Gamma_{D,0}^u \cup \Gamma_{N,0}^u}=
\overline{\Gamma_{D,0}^t \cup \Gamma_{N,0}^t}= \partial\Omega_0$. 
The weak forms of the balance laws in~\eqref{eqn:LMB_AMB} 
and~\eqref{eqn:energy_bal} are expressed respectively as 
\begin{align}\label{eqn:weakform}
\begin{split} 
    \int_{\Omega_0} \xibold\cdot\rho_0\abm\,dV + 
    \int_{\Omega_0} \bnabla^s\xibold\cdot\sigmabold \,dV\ &=\ 
\int_{\Gamma_{N,0}^u} \xibold\cdot\bar\pbm \,dA
    +\int_{\Omega_0} \xibold\cdot\rho_0\bbm \,dV\ , \\
    \int_{\Omega_0} \zeta c \dot\theta\,dV - 
\int_{\Omega_0} \nabla\zeta\cdot\qbm_0 \,dV\ & =\ 
    \int_{\Omega_0} \zeta \rho_0 r \,dV + \int_{\Omega_0} \zeta\theta
\Mbm\compactcdot\dot\epsbold \,dV - \int_{\Gamma_{N,0}^t} \zeta\bar q \,dA\ , 
\end{split} 
\end{align}
where~$\bnabla^s$ denotes the symmetric part of the referential 
gradient, while~$\bar\pbm$ and~$\bar q$ are the imposed traction and 
flux boundary conditions on~$\Gamma_{N,0}^u$ and~$\Gamma_{N,0}^t$, respectively.
The weak forms are expressed in terms of arbitrary and sufficiently smooth 
weighting functions~$\xibold$ for the linear momentum balance and~$\zeta$ 
for the energy balance, each satisfying the corresponding homogeneous
Dirichlet boundary conditions. The infinitesimal strain tensor is related
to the displacement field~$\ubm$ according to~$\epsbold=\bnabla^s\ubm$. 
\par
Taking into account the arbitrariness of the weighting functions and
introducing standard displacement-type finite element piecewise approximations 
for the dependent variables, their derivatives, and the corresponding 
weighting functions leads to elemental equations written in matrix form as 
\begin{equation}\label{eqn:discrete-balance-ele}
 \centering
 \begin{aligned}
   \left[\Mbm^e_{u}\right]\left[\hat{\abm}^e_{n}\right] +\left[ \Rbm_{u,n}^e\right] 
   & \ =\ \left[\Fbm_{n}^e\right]+
\int_{\partial\Omega_0^e \setminus
\Gamma_{N,0}^u}\left[\Nbm_u^{e}\right]\left[\pbm_{n}^e\right]\,dA\ ,\\
   \left[\Tbm^e\right]\left[\hat{\dot{\thetabold}}_{n}^e\right] + \left[\Mbm^e_{t}\right]\left[\hat{\thetabold}_{n}^e\right]+ 
    \left[\Rbm^e_{t,n} \right]
   & \ =\  
   \left[\Qbm_{n}^e\right]-
\int_{\partial\Omega_0^e\setminus\Gamma_{N,0}^t}
\left[\Nbm_t^{e}\right]^T \left[\qbm_{n}^e\right]\,dA\ .
 \end{aligned}
\end{equation}
Here, the subscript ``$n$'' identifies an algebraic quantity as estimated 
at time~$t_{n}$, while the superscript~``$e$'' refers to the element~$e$. 
In addition, the overset symbol in~$\hat{(\cdot)}$ implies that the term under 
it comprises the set of nodal values associated with the quantity~$(\cdot)$. 
All the matrices used in~\eqref{eqn:discrete-balance-ele} are defined 
in the appendix. Upon neglecting, as is customary, the interelement flux terms, 
the preceding elemental equations give rise to the global system, expressed as 
\begin{equation}\label{eqn:discrete-balance-global}
 \centering
 \begin{aligned}
 [\Mbm_{u}] \left[\hat{\abm}_{n}\right] +\left[ \Rbm_{u,n}\right]
 - [\Fbm_{n}]& \ =\ [\zerobold]\ ,\\
    \left[\Tbm\right]\left[\hat{\dot{\thetabold}}_{n}\right] + \left[\Mbm_{t}\right]\left[\hat{\thetabold}_{n}\right]+ 
 \left[\Rbm_{t,n}\right]- \left[\Qbm_{n}\right]& \ =\ [\zerobold]\ , 
 \end{aligned}
\end{equation}
where any elemental quantity~$[(\cdot)^e]$ is assembled into its
global counterpart~$[(\cdot)]$ by means of a standard assembly operation. 
\par
The temporal discretization of the displacement vector~$[\hat{\ubm}]$ and 
the velocity vector~$[\hat{\vbm}]$ is effected 
by the standard trapezoidal rule, while the backward Euler rule
is employed for the temperature vector~$[\hat{\thetabold}]$. 
Therefore, given a time interval~$(t_{n-1},t_{n}]$ of 
size~$\Delta t_{n-1} = t_{n}-t_{n-1}$,
\begin{equation} \label{eqn:time-discretization}
 \centering
  \begin{aligned}
    \left[{\hat\ubm}_{n}\right]\ &=\ \left[\hat{\ubm}_{n-1}\right] + 
\left[\hat{\vbm}_{n-1}\right]\dt{n-1} + \frac{1}{4}
    \left\{\left[\hat{\abm}_{n-1}\right]+\left[\hat{\abm}_{n}\right]\right\}\dt{n-1}^2\ ,\\
    \left[\hat{\vbm}_{n}\right]\ &=\ \left[\hat\vbm_{n-1}\right] + 
\frac{1}{2}\left\{\left[\hat\abm_{n-1}\right]+\left[\hat{\abm}_{n}\right]\right\}\dt{n-1}\ ,\\
    \left[\hat\thetabold_{n}\right]\ &=\
\left[\hat\thetabold_{n-1}\right] +
\left[\hat{\dot\thetabold}_{n}\right]{\dt{n-1}}\ .
  \end{aligned}
\end{equation}
\par
Taking into account the constitutive 
laws~(\ref{eqn:fourier_law}--\ref{const-elastic-stress}) and the time-stepping 
rules in~\eqref{eqn:time-discretization}, the system of linear
algebraic equations at the~$k$-th Newton-Raphson iteration in the time 
interval~$(t_{n-1},t_{n}]$ can be written in matrix form as 
\begin{equation} \label{eqn:disc-matrix-Newton}
 \centering
  \begin{aligned}
   &\left[\Mbm_{u}\right]\left[\hat{\abm}^{(k)}_{n}\right]
+ \left[\Rbm_{u,n}^{(k)}\right] - \left[\Fbm_{n}\right] \\
   & \qquad\qquad\qquad\qquad\qquad+\left[\Kbm^{(k)}_{u,n}\right]
\left[\Delta\hat{\ubm}^{(k)}_{n}\right] +
\left[\Kbm^{(k)}_{t,n}\right] 
\left[\Delta \hat{\thetabold}^{(k)}_{n}\right]\ =\ [\zerobold]\ ,\\
   &\left[\Tbm\right]\left[\hat{\dot{\thetabold}}_{n}^{(k)}\right] + 
\left[\Mbm_{t}\right]\left[\hat{\thetabold}_{n}^{(k)}\right]+ 
    \left[\Rbm^{(k)}_{t,n}\right]
    - \left[\Qbm_{n}\right]  \\ 
   & \qquad\qquad\qquad\qquad\qquad+\left[\Abm^{(k)}_{u,n}\right]\left[\Delta\hat{\ubm}_{n}^{(k)}\right]+\left[\Abm^{(k)}_{t,n}\right]\left[\Delta
\hat{\thetabold}_{n}^{(k)}\right]\ =\ [\zerobold]\ ,
  \end{aligned}
\end{equation}
where~$\Delta(\cdot)_{n}$ is the change of the variable~$(\cdot)$ 
in~$(t_{n-1},t_{n}]$. Expressions for the tangent 
matrices~$\left[\Kbm^{(k)}_{u,n}\right]$, $\left[\Kbm^{(k)}_{t,n}\right]$, 
$\left[\Abm^{(k)}_{u,n}\right]$, and~$\left[\Abm^{(k)}_{t,n}\right]$ 
in~\eqref{eqn:disc-matrix-Newton} are provided in the appendix. 
The global displacement and temperature vectors are updated according to 
\begin{equation}\label{eq:tint}  
[\hat\ubold_{n}^{(k+1)}]\ =\  
[\hat\ubold_{n}^{(k)}] + 
[\Delta\hat\ubold_{n}^{(k)}]\quad,\quad 
[\hat\thetabold_{n}^{(k+1)}]\ =\  
[\hat\thetabold_{n}^{(k)}] + 
[\Delta\hat\thetabold_{n}^{(k)}]\ , 
\end{equation} 
until convergence is attained in a predefined stopping criterion and 
a prescribed tolerance. 
\par
To accommodate new nodes and elements generated when material is deposited 
on the boundary of the existing body, each node is assigned an initial
displacement and temperature. These ``history'' variables ensure zero 
initial stress and strain. To this end, the displacement history
values are subtracted from the nodal displacement values when
computing the strain and stress at a typical integration point. 
Likewise, the temperature history value represents the reference 
temperature~$\theta_0$ in~\eqref{const-elastic-stress} of an element when 
it is first created, computed as the average of the nodal temperature of the element. 

For printed objects that have curved boundaries, our simulation relies on a slicing algorithm 
to determine how the curves are approximated and printed. The choice of the slicing algorithm is not a focus of this study. 
For the one that is used in our finite element model, the reader can find it in~\cite{wang2021}. 
In general however, curved boundaries require the use of partial elements that involve hanging nodes.   
The reference coordinates of the hanging nodes are determined by the slicing algorithm while their displacement 
and temperature are solved together with the rest of the displacement and temperature field. 
To do so, additional conditions need to be imposed if a hanging node lies on an existing edge.
When such a node is created, a ratio of its coordinates and those of two nodes that define this existing edge can be computed. 
This ratio is then maintained when computing for the displacement and temperature of the hanging node, so that it would remain at the same relative position of the edge throughout the simulation. 
Full details on the implementation of history variables and hanging nodes 
can be found in~\cite{wang2021}. 
\section{Optimization problem formulation}\label{se:opt}
\par\noindent
Let~$f(\ubm,\thetabold;\ybm,t)$ be the objective function of the design 
optimization problem, where 
$\ybm\in{\Bbb R}^d$ represents the vector set of~$d$ design variables,  
\textit{i.e.}, process parameters such as chamber temperature and 
layer thickness. The continuum-level optimization problem can be expressed as 
\begin{equation} \label{eqn:opt-fcn}
 \centering
  \begin{aligned}
   &\underset{\substack{\ybm\in{\Bbb R}^d}}{\text{minimize}} 
   & &f\left(\ubm,\thetabold;\ybm,t\right)\ ,\\
   &\text{subject to}
   & &\rbm\left( \ubm,\thetabold;\ybm,t\right)\ =\
\zerobold\ ,\quad c_i(\ybm)\ \leqslant\ 0\ ,\ i = 1, \ldots, m\ . 
  \end{aligned}
\end{equation}
Here, $\rbm$ comprises the balance laws in~\eqref{eqn:weakform}
(before the finite element approximation), while~$c_i(\cdot)$, $i=1,\ldots,m$, 
are inequality constraint functions on the design variables. 
\par
\begin{figure}
 \centering
  \includegraphics{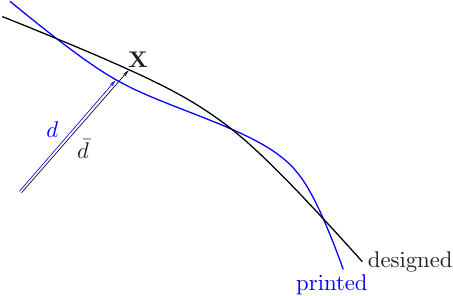}
 \caption{Schematic of the distance functions~$d$ and~$\bar{d}$} 
\label{fig:dd}
\end{figure}

In this paper, the objective function is chosen to express the deviation of
the part's actual shape from the desired shape. Such deviation, or shape error, may pose 
significant challenges in AM~\cite{pandey2003}. In the two-dimensional case, 
given a fixed reference configuration, shape error may be defined as 
\begin{equation} \label{eqn:obj-funcl}
   f\left(\ubm,\thetabold;\ybm,t\right) \ =\ 
\frac{\displaystyle\left[\int_{\Gamma_{S}^{t}} 
\left(d\left(\Xbm,\ubm\right)-\bar d(\Xbm)\right)^2 dl\right]^{\frac{1}{2}}}
          {\displaystyle L_c \left(\int_{\Gamma_{S}^{t}} ~dl\right)^{\frac{1}{2}}}\ ,
\end{equation}
where~$\Xbm$ represent the placement of the boundary points on
the (necessarily time-dependent) surface of interest~$\Gamma_S^t$
in the designed geometry under no deformation. Also, 
${d\left(\Xbm,\ubm\right)-\bar{d}(\Xbm)}$ is a problem-dependent scalar 
measure of the difference of the printed position~$d$ from the designed
position~$\bar{d}$ on~$\Gamma_S^t$, while~$L_c$ is a representative length 
of the part intended to render~$f$ dimensionless.   
As seen from~\eqref{eqn:obj-funcl}, the objective function depends explicitly 
on~$\ubm$ alone. Upon finite element discretization, the objective function 
$f$ is measured on one-dimensional element edges 
and can be evaluated numerically using Gaussian quadrature. The discrete counterpart~$f_h$ of the objective function for linear elements 
is given by 
\begin{equation} \label{eqn:obj-funcl-disc}
 \centering
  \begin{aligned}
  f_h(\hat{\ubm},\hat{\thetabold};\by,t)\ =\ 
\frac{\displaystyle\left[ \sum_{i=1}^p \sum_{j=1}^l \left(d\left(\Xbm_j,\ubm_j
\right)-\bar d(\Xbm_j)\right)^2 w_j \frac{\Delta l_i}{2}\right]^
            {\frac{1}{2}}}
            {\displaystyle L_c \left(\sum_{i=1}^p \Delta l_i\right)^{\frac{1}{2}}}\ ,
  \end{aligned}
\end{equation}
where~$p$ is the number of discretized edges of interest, $\Delta l_i$ is 
the length of the~$i$-th edge, $l$ is the number of quadrature points per edge, 
and~$w_j$ is the weight of the~$j$-th quadrature point. 
The values of $\Xbm_j,\ubm_j$ at quadrature points are evaluated using $\Xbm$ and the nodal displacement $\hat{\ubm}$.
The precise definition 
of the position deviation~$d(\ubm)-\bar{d}$ is specified 
in Section~\ref{se:exp} for each example. Implicit in the 
preceding definitions of~$f$ and~$f_h$ is that~$d$ induces a bijection 
between the designed to the printed surface regardless of the value of the
displacement~$\ubm$, see Figure~\ref{fig:dd}.
\par
Depending on the nature of $\ybm$, the optimization 
problem~\eqref{eqn:opt-fcn} may require the use of substantially different 
algorithms. For example, some process parameters, such as chamber temperature,
vary continuously. Hence, the objective function~$f_h$ may be considered 
implicitly differentiable in them. This would enable the application 
of well-established gradient-based optimization algorithms. 
On the other hand, parameters such as layer thickness, which is mandated 
by the (integer) number of layers for a given specimen size, are discrete and, 
therefore, non-differentiable. For optimization in such variables, 
a gradient-free algorithm would be appropriate. 
%
%
\section{Gradient-descent method with line search}\label{se:grad}
\par\noindent
For gradient-based optimization methods, the derivatives of the 
discretized objective~$f_h$ from~\eqref{eqn:obj-funcl-disc} with respect to 
design variables~$\ybm$ must be known or estimated. Let~$N$ be the total number 
of discrete time steps in the analysis. For simplicity, denote the
combined system of fully discretized PDEs 
at time~$t_{n},1\leq n\leq N$, as~$\rbm^h_n$, where the superscript~$h$ 
is included to emphasize the discrete vector form of the variables. 
It follows that the discretized form of the PDE constraint 
in~\eqref{eqn:opt-fcn} can be expressed as 
\begin{equation}\label{eqn:pde-general}
 \centering
\rbm^h_n \left( \hat\ubm_{n}^h, \hat\ubm_{n-1}^h, \hat\vbm_{n-1}^h, 
\hat\thetabold_{n}^h, \hat\thetabold_{n-1}^h, 
\hat\ubm^{his}_{n}, \hat\thetabold^{his}_{n};\ybm,t_{n}\right)\ =\ \zerobold\ ,
\end{equation}
in view of the single-step time integration rules in~\eqref{eq:tint}. 
As noted in Section~\ref{se:goveqn}, the history variables 
$\ubm^{his}_{n},\thetabold^{his}_{n}$ represent the displacement and temperature field at each new node and element 
when they are added at their respective 
time step. Thus, they are also discrete in form and full details on their 
precise definition may be found in~\cite{wang2021}. 
In addition, it is assumed that any inequality constraints~$c_i(\cdot)$ 
on~$\ybm$ can be included by requiring that~$\ybm \in \mathfrak{C}$, 
where~$\mathfrak{C}\subset \Rbb^m$ is a suitably defined convex and compact set.
\par
The solver-consistent discretized form of~\eqref{eqn:opt-fcn} is now expressed 
as 
\begin{equation} \label{eqn:opt-fcn-reduced}
 \centering
  \begin{aligned}
  &\underset{\substack{\ybm \in \mathfrak{C}}}{\text{minimize}} 
  & &f_h\left(\hat\ubm^h_{N},\hat\thetabold^h_{N};\ybm,t_N\right)\\
   &\text{subject to}
  & &\rbm^h_n \left(
\hat\ubm^h_{n},\hat\ubm^h_{n-1},\hat\vbm_{n-1}^h,
\hat\thetabold^h_{n},\hat\thetabold^h_{n-1},
\hat\ubm_{n}^{his},\hat\thetabold_{n}^{his};\ybm,t_n\right)\
=\ \zerobold\quad ,\quad n=1,2,\ldots,N\ .
  \end{aligned}
\end{equation}
The preceding discrete optimization problem reflects the fact that
only the end-state is of consequence in the design.
\par
Consider the reduced-space formulation, in which~$\ubm^h_N$ and 
$\thetabold^h_N$ are regarded as implicit functions of~$\ybm$ 
at time~$t_N$ and also let~$N$ be independent of~$\ybm$. 
Then, by chain rule, the gradient of~$f_h$ with respect to~$\ybm$ is 
\begin{equation}\label{eqn:obj-diff}
 \centering
 \centering
  \frac{d f_h}{d \ybm} \ =\  
\pder{f_h}{\hat\ubm_N^h}\pder{\hat\ubm_N^h}{\ybm} +
\pder{f_h}{\hat\thetabold_N^h}\pder{\hat\thetabold_N^h}{\ybm} + 
\pder{f_h}{\ybm}\ .
\end{equation}
Each term on the right-hand-side of~\eqref{eqn:obj-diff} is derived next. 
First, note that, owing to the explicit dependence of~$f_h$ 
on~$\ybm$, $\hat\ubm_N^h$, and $\hat\thetabold_N^h$, the terms  
$\pder{f_h}{\ybm}$, $\pder{f_h}{\hat\ubm_N^h}$, and 
$\pder{f_h}{\hat\thetabold_N^h}$ can be computed directly by differentiating 
the function~$f_h$. Of course, given the specific form of
$f_h$ in~\eqref{eqn:obj-funcl-disc}, the only non-zero such term here is
$\pder{f_h}{\hat\ubm_N^h}$. On the other hand, the state sensitivities 
$\pder{\hat\ubm_N^h}{\ybm}$ and $\pder{\hat\thetabold^h_N}{\ybm}$ at~$t_N$ 
are dependent on those of the previous time steps, that is, 
$\pder{\hat\ubm_n^h}{\ybm}$ and $\pder{\hat\thetabold^h_n}{\ybm}$, 
for $n=1,\ldots,N-1$. Due to the small total number of optimization variables, 
a pure primal approach is adopted here to obtain the sensitivities.
To this end, upon differentiation of both sides of 
Equation~\eqref{eqn:pde-general} with respect to~$\ybm$, it follows that 
\begin{multline}\label{eqn:pde-diff}
   \pder{\rbm_{n}^h}{\hat\ubm_{n}^h}\pder{\hat\ubm_{n}^h}{\ybm}+\pder{\rbm_{n}^h}{\hat\ubm_{n-1}^h}\pder{\hat\ubm_{n-1}^h}{\ybm}+
   \pder{\rbm_{n}^h}{\hat\vbm_{n-1}^h}\pder{\hat\vbm_{n-1}^h}{\ybm}+\pder{\rbm_{n}^h}{\hat\thetabold_{n}^h}\pder{\hat\thetabold_{n}^h}{\ybm}+
   \pder{\rbm_{n}^h}{\hat\thetabold_{n-1}^h}\pder{\hat\thetabold_{n-1}^h}{\ybm} \\
   +\pder{\rbm_{n}^h}{\hat\ubm^{his}_{n}}\pder{\hat\ubm^{his}_{n}}{\ybm} +
   \pder{\rbm_{n}^h}{\hat\thetabold^{his}_{n}}\pder{\hat\thetabold^{his}_{n}}{\ybm}+
   \pder{\rbm_{n}^h}{\ybm} \ =\ \zerobold\ .  
\end{multline}
The state sensitivities $\pder{\hat\ubm_{n}^h}{\ybm}$ and
$\pder{\hat\thetabold_{n}^h}{\ybm}$ at time~$t_n$ are obtained 
from~\eqref{eqn:pde-diff}, in terms of the rest of the (known) terms 
in the equation, including the state sensitivities at~$t_{n-1}$. 
The initial conditions for these sensitivities are set to zero, since the
initial displacement and temperature fields are prescribed. 
\par
It is important to note here that use of the Newton-Raphson method
for the solution to $\rbm_n^h=\zerobold$ requires the calculation of the
``stiffness'' terms $\pder{\rbm_{n}^h}{\hat\ubm_{n}^h}$ and 
$\pder{\rbm_{n}^h}{\hat\thetabold_{n}^h}$. Therefore, these
two terms (which are derived in the appendix) are already available 
for the sensitivity calculations. Furthermore, the terms 
$\pder{\rbm_{n}^h}{\hat\ubm_{n-1}^h}$, $\pder{\rbm_{n}^h}{\hat\vbm_{n-1}^h}$, 
and $\pder{\rbm_{n}^h}{\hat\thetabold_{n-1}^h}$, which appear
in~\eqref{eqn:pde-diff}, can be readily obtained by substituting 
the time-stepping formulae~\eqref{eqn:time-discretization} 
into~\eqref{eqn:discrete-balance-global} and taking the derivatives
of the discrete governing equations comprising $\rbm_n^h$ with respect 
to~$\ubm_{n-1}^h,\thetabold_{n-1}^h,\vbm_{n-1}^h$. 
Again, explicit expressions for these terms may be found in the appendix.
\par
Based on the implementation of the history variables $\hat\ubm_{n}^{his}$
and $\hat\thetabold_n^{his}$, it is straightforward to compute the
corresponding sensitivities $\pder{\hat\ubm^{his}_{n}}{\ybm}$ and
$\pder{\hat\thetabold^{his}_{n}}{\ybm}$ in Equation~\eqref{eqn:pde-diff} 
at the element level and then assemble into their global counterparts 
(for details see, again,~\cite{wang2021}). 
Since the history variables are simply the displacement and temperature 
of a node and an element, respectively, at the time they are created, they remain constant throughout 
the simulation.  Thus, the displacement history variable sensitivity 
$\pder{\hat\ubm^{his,e}_{n}}{\ybm}$ is simply the displacement sensitivity at the time of the creation of a node,
and the temperature history sensitivity $\pder{\hat\thetabold^{his,e}_{n}}{\ybm}$ is the average 
of the nodal temperature sensitivity of an element. 
These sensitivities are stored for subsequent use. The same applies to the derivative
terms $\pder{\rbm_{n}^h}{\hat\ubm^{his}_{n}}$ and 
$\pder{\rbm_{n}^h}{\hat\thetabold^{his}_{n}}$, which are derived in the
appendix.
\par
At each time~$t_n$, the state sensitivities at $t_{n-1}$ 
are stored and used in the solution of~\eqref{eqn:pde-diff}.  
%
The last remaining term in~\eqref{eqn:pde-diff}, $\pder{\rbm_{n}^h}{\ybm}$,
can be explicitly computed once the optimization variable $\ybm$ is defined. 
With all the necessary terms now available, Equation~\eqref{eqn:pde-diff} 
is solved for $\pder{\hat\ubm_n^h}{\ybm}$ and $\pder{\hat\thetabold^h_n}{\ybm}$,
subject to the initial sensitivities (taken here to be zero). 
\par
All the sensitivities computed here analytically have been tested against 
second-order accurate centered-finite difference approximation to ensure 
the faultlessness of the derivations.
%
\par
Given the initial displacement vector~$\hat\ubm_0^h$, temperature 
vector~$\hat\thetabold_0^h$ and design variable vector~$\ybm_0$, 
it is now possible to compute the gradient~\eqref{eqn:obj-diff}. 
The optimization problem~\eqref{eqn:opt-fcn-reduced} can be recast 
in reduced-space as 
\begin{equation} \label{eqn:opt-fcn-2}
 \centering
  \begin{aligned}
   &\underset{\substack{\ybm \in \mathfrak{C}}}{\text{minimize}} 
   & &f_h\left(\ubm_0,\thetabold_0,\ybm,t_N\right)\ .
  \end{aligned}
\end{equation}
For relatively small dimension of~$\ybm$, many well-established 
nonlinear optimization methods may be used to solve this problem. 
Here, a projected gradient-descent method with line search is employed, 
as described in Algorithm~\ref{alg:descent}. The initial sensitivities 
$\pder{\hat\ubm_0}{\ybm}$ and $\pder{\hat\thetabold_0}{\ybm}$ are set 
to~$0$ as they do not depend on~$\ybm$ at the starting time~$t_0$.
\begin{algorithm}
 \DontPrintSemicolon
 \SetAlgoNoLine 
 \caption{Projected gradient descent with line search}\label{alg:descent}
  Let~$\alpha_0=1$, $\rho\in(0,1)$, $\eta\in(0,1)$\; 
  \For{$k=0,1,2,...$}{
     Perform the finite element simulation. Compute sensitivities and the gradient in~\eqref{eqn:obj-diff}.\;
     Let $\pbm_k = -\frac{d f_h\left(\ybm_k\right)}{d \ybm_k}/ \| \frac{d f_h\left(\ybm_k\right)}{d \ybm_k} \| $. \;
     Set $\alpha_k \leftarrow \alpha_0$.\;
     \While{$f_h \left(\ybm_k + \alpha_k \pbm_k \right) > f_h\left(\ybm_k \right)+
            \eta \alpha_k \left[\frac{d f_h\left(\ybm_k\right)}{d \ybm_k}\right]^T \pbm_k$ }{
       $\alpha_k \leftarrow \rho\alpha_k $\;
     }
     $\ybm_{k+1} = \ybm_k + \alpha_k \pbm_k$\;
     Project $\ybm_{k+1}$ onto $\mathfrak{C}$.
  }
\end{algorithm}
\par
In this algorithm, $\eta$ and $\rho$ are user-specified constants,
while the scalar~$\alpha_k$ is the length of the line-search step. 
It is worth pointing out that the line-search 
step could be computationally expensive as~$f_h$ needs to be evaluated 
each time~$\alpha_k$ is updated. The projection step for a general convex 
and compact~$\mathfrak{C}$ typically involves solving 
a minimum-distance problem. However, in the case of bound constraints, 
one is only required to force $\ybm$ to stay inside the bounds.
For more complicated constraints which are nonlinear in~$\ybm$, this algorithm 
would not suffice. For the purposes of this work, Algorithm~\ref{alg:descent} 
is implemented in both C++ and Python. 
%
%
\section{Gradient-free Optimization Methods}\label{se:grad-free}
\par\noindent
Two optimization algorithms are presented in this section for the solution 
of~\eqref{eqn:opt-fcn} and~\eqref{eqn:opt-fcn-2} that do not require 
sensitivities and rely only on objective function evaluations. 
The first is a local search algorithm that advances by comparing 
objective values at neighboring points in the design space. 
Owing to parallel computing, this algorithm can quickly yield reasonable 
results and is easy to implement. The second is a Bayesian optimization
algorithm, where a surrogate model for the objective is built using samples. 
The solution to the optimization problem is updated based on the prediction 
of the surrogate model. Algorithms of this type have received significant 
attention in the past decade and have been shown to efficiently generate 
global optimal result~\cite{frazier2018tutorial}. 
\subsection{Method of local variations}\label{sec:loc-var}
This is a version of the method of local variations proposed by 
Polak~\cite{polak} and is shown in Algorithm~\ref{alg:lv}. 
It is used here to solve the reduced formulation~\eqref{eqn:opt-fcn-2} of 
the design optimization problem. The algorithm starts with a guess 
in the optimization variable space and takes a series of steps towards 
an optimal solution by constantly comparing neighborhood objective values 
and adjusting step sizes. Multiple starting points can be used 
to improve the results, as the algorithm can be shown to reach 
a local minimum if the objective is continuously differentiable~\cite{polak}. 
\begin{algorithm}
 \DontPrintSemicolon
 \SetAlgoNoLine 
   \caption{Method of local variations }\label{alg:lv}
     Select a $\ybm_0 \in \Rbb^d$ such that $\ybm_0 \in \mathfrak{C}$ is
feasible. Select a $\boldsymbol{\tau}_0 > 0$.\;
     Set $k=0$, and compute $f_h(\ybm_0)$. Let $f_{min} = f_h(\ybm_0)$.;
     \While{$\boldsymbol{\tau}_{k} > \boldsymbol{\tau}_{min}$}{
       Set $f_{min} = f_h(\ybm_k)$ and $\dbm_{min}=0$.\;
      \For{$j =  1,\dots, 2n$}{
         Check if $\zbm_k^j = \ybm_k + (\boldsymbol{\tau}_k)_{j/2} \dbm_j \in \mathfrak{C}$. Project
it onto $\mathfrak{C}$ if not so that $\zbm_k^j\in \mathfrak{C}$.\; 
         Compute $f_h( \zbm_k^j)$.\; 
         \If{$f_h\left( \zbm_k^j \right) < f_{min}$}{
           $f_{min} = f_h ( \zbm_k^j )$. \;
           $\dbm_{min} = \zbm_k^j - \ybm_k$.\;
         }
      }
      $\ybm_{k+1} = \ybm_k + \dbm_{min}$.\;
      \If{$f_{min}<f_h(\ybm_k)$}{
          $\boldsymbol{\tau}_{k+1} = \boldsymbol{\tau}_k$, $k = k+1$.\;
      }
      \Else{
         Reduce the step size via $(\boldsymbol{\tau}_{k+1})_j = (\boldsymbol{\tau}_k)_j /2$ or $(\boldsymbol{\tau}_{k+1})_j = (\boldsymbol{\tau}_k)_j - 1$ for $j=1,\dots,n$. Set $k=k+1$.\;
      }
     }
\end{algorithm}
Denote the initial step size as $\boldsymbol{\tau}_0$  and start with $\ybm_0\in \mathfrak{C}\subset \Rbb^d$. 
The minimum step size is controlled by vector $\boldsymbol{\tau}_{min}$. 
At each iteration~$k$, the algorithm centers around the current point~$\ybm_k$ 
in the optimization space and runs $2d$ simulations on neighboring points, with the neighboring points 
selected by the step size $\boldsymbol{\tau}$ and 
directions $\dbm_j$, $j=1,\cdots,2d$. If neighboring points
produce smaller values of the objective, the point that produces the smallest 
value is selected to be the next iteration point. 
If no such point can be found, the current step size
$\boldsymbol{\tau_k}$ is reduced and simulation at new neighboring points are run until
the step size of each variable becomes smaller than each component of $\boldsymbol{\tau}_{min}$. 
The quantity $\dbm_j$ represents the increasing and decreasing search direction 
for each component of~$\ybm$. In Algorithm~\ref{alg:lv}, the default 
optimization variable is continuous, and thus the step size can be repeatedly 
reduced by half. However, the algorithm can be readily adapted 
to discrete optimization variables. Specifically, for an integer parameter 
such as the number of deposition layers, one may reduce the 
step size $\norm{\boldsymbol{\tau}_k \dbm_j}$ by 1. 
The algorithm enforces the simple constraints at step $6$ so that all points evaluated are feasible. 
\subsection{Bayesian optimization method}\label{se:bayesian}
\par\noindent
A Bayesian optimization algorithm~\cite{frazier2018} is alternatively employed 
to solve Problem~\eqref{eqn:opt-fcn-2}. 
This is a ``black-box'' sequential surrogate-based optimization method, 
where the objective and the associated PDEs are represented by 
a surrogate model and, therefore, no sensitivity calculations are necessary. 
The surrogate model is typically a random function with a prior that 
is constructed and updated based on the samples (data points) collected. Then, a posterior distribution is formed using 
the sequentially updated prior so that a prediction of the objective 
at a given point can be made~\cite{frazier2018tutorial,shahriari2015taking}. 
In addition, an acquisition function is used to determine where to sample next. 
This function strikes a balance between exploitation and exploration, 
which correspond to an improved objective value and lower uncertainty, 
respectively. 
\par
A sample in the Bayesian optimization process is the individual 
finite element simulation and its objective value given the process 
parameters~$\ybm$. The well-tested Gaussian process (GP) is chosen 
for the surrogate model, while the expected improvement acquisition 
function is adopted. Given the sample points, 
GP assumes a multivariate jointly-Gaussian distribution between~$\ybm$ 
and the objective function~$f_h$: 
\begin{equation} \label{eqn:GP-1}
 \centering
  \begin{aligned}
  \begin{bmatrix}
   f_h(\ybm_1) \\
   \vdots\\
           f_h(\ybm_N)
  \end{bmatrix} \sim \mathcal{N}\left( 
      \begin{bmatrix}    
                 m(\ybm_1)\\
 \vdots\\
 m(\ybm_N)
      \end{bmatrix},
              \begin{bmatrix}
               k(\ybm_1,\ybm_1) \dots k(\ybm_1,\ybm_N)\\
      \vdots\\
               k(\ybm_N,\ybm_1) \dots k(\ybm_N,\ybm_N)
      \end{bmatrix}
      \right)\ . 
\end{aligned}
\end{equation}
Here, $\mathcal{N}$ is the normal distribution, $\ybm_1,\ldots,\ybm_N$ are the samples, $m(\cdot):\Rbb^N\to\Rbb$ is the mean function, 
and~$k(\cdot,\cdot):\Rbb^N\times\Rbb^N\to\Rbb$ is the covariance function. Then, the
posterior probability distribution of a new~$\ybm$ can be inferred
using Bayes' rule~\cite{frazier2018},
\begin{equation} \label{eqn:GP-post}
 \centering
  \begin{aligned}
  &f_h(\ybm) | \ybm, \ybm_{1:N},f_h(\ybm_{1:N}) \sim \mathcal{N}
(\mu(\ybm),\sigma^2(\ybm))\ ,\\
  &\mu(\ybm)\ =\ k(\ybm,\ybm_{1:N}) k(\ybm_{1:N},\ybm_{1:N})^{-1}
\left(f_h(\ybm_{1:N})-m(\ybm_{1:N}) \right) + m(\ybm_{1:N})\ ,\\
  &\sigma^2(\ybm)\ =\
k(\ybm,\ybm)-k(\ybm,\ybm_{1:N})k(\ybm_{1:N},\ybm_{1:N})^{-1}\sigma(\ybm_{1:N},\ybm)\ ,
\end{aligned}
\end{equation}
where the vector~$\ybm_{1:N}$ is the notation for $\ybm_1,\dots,\ybm_N$ and 
$k(\ybm_{1:N},\ybm_{1:N}) = [k(\ybm_1,\dots,k(\ybm_N) ;$ $\dots; k(\ybm_N,\ybm_1),\dots,k(\ybm_N,\ybm_N)]$. 
The function $\mu(\cdot)$ and $\sigma^2(\cdot)$ are referred to as the posterior mean and variance, respectively.
The mean function $m(\cdot)$ is often simply a constant function.
The covariance functions, or kernels, of the GP is critical to an accurate 
surrogate model. The Squared Exponential Covariance Function~\cite{brochu2010} 
or power exponential kernel is adopted in this paper. The definition of the 
covariance function is 
\begin{equation} \label{eqn:kernel-sec}
 \centering
  \begin{aligned}
  k(\ybm,\ybm';\theta)\ =\
\text{exp}\left(-\frac{\norm{\ybm-\ybm'}^2}{\theta^2}\right)\ ,
  \end{aligned}
\end{equation}
where~$\theta$ denotes the hyper-parameters of the kernel.
The GP model's hyper-parameters are optimized during the regression 
by maximizing the log-marginal-likelihood with 
the BFGS method~\cite{optimization}. The expected improvement acquisition function 
can be written as
\begin{equation} \label{eqn:EI-1}
 \centering
  \begin{aligned}
       EI(\ybm)\ =\ \begin{cases}
       (\mu(\ybm) - f_h(\ybm^+) - \xi)\Phi(Z)+\sigma(\ybm)\phi(Z), &\text{if} \ \sigma(\ybm)>0\\
       0, &\text{if} \ \sigma(\ybm) = 0    
                 \end{cases}\ ,
  \end{aligned}
\end{equation}
where 
\begin{equation} \label{eqn:EI-2}
 \centering
  \begin{aligned}
       Z\ =\ \begin{cases}
       \frac{\mu(\ybm)-f_h(\ybm^+)-\xi}{\sigma(\ybm)}, &\text{if} \ \sigma(\ybm)>0\\
       0, &\text{if} \ \sigma(\ybm) = 0\ 
                 \end{cases},
  \end{aligned}
\end{equation}
and $\ybm^+=\text{argmax}_{\ybm_i\in\ybm_{1:N}}f_h(\ybm_i)$ of all~$N$ current 
samples~\cite{brochu2010}. The trade-off parameter~$\xi>0$ controls the balance between exploration and exploitation. At an input~$\ybm$, the mean value predicted by 
the surrogate is~$\mu(\ybm)$ and the variance is~$\sigma(\ybm)$. The functions 
$\phi$ and~$\Phi$ are the probability density function (PDF) and cumulative 
distribution function (CDF) of the standard normal distribution, respectively. 
In order to find the maximal point of the acquisition function, which would be 
the next sample point, optimization algorithms including L-BFGS or random 
search can be used. The Bayesian optimization (BO) algorithm applied 
to~\eqref{eqn:opt-fcn-2} is given in Algorithm~\ref{alg:bo}. 
An example convergence measure that allows the algorithm to terminate is if consecutive optimization iterations produce
optimal point and optimal value within a tolerance range, using the updated Gaussian process
surrogate model that includes the latest samples.
\begin{algorithm}
 \DontPrintSemicolon
 \SetAlgoNoLine 
 \caption{Bayesian optimization with Gaussian process}\label{alg:bo}
  Choose initial sampling data $\ybm_0$.\; Train the Gaussian process model with $\ybm_0$.\; 
  \For{$k=1,2,\dots$}{
     Evaluate the acquisition function using sample points $\{\ybm_0,\dots,\ybm_{k-1}\}$ and find $\ybm_k = \text{argmax}_{\ybm} EI(\ybm)$.\;
     Run finite element simulation with design variables $\ybm_k$. \;
     Compute the objective $f_h(\ybm_k)$ based on the simulation result. \;
     Retrain the GP surrogate model with the addition of the new sample $\ybm_k$ and $f_h(\ybm_k)$.\;
     Solve the optimization problem  with the updated surrogate model. \;
     Evaluate convergence measure. Exit if satisfied. \;
  }
\end{algorithm}
\par
Algorithm~\ref{alg:bo} is implemented in Python and C++ with scikit-learn 
machine learning library~\cite{scikit-learn,sklearn_api} and uses 
standardscaler preprocessing which reduces the mean of the objective data to $0$ and scales them to unit variance. We note that Bayesian optimization is typically implemented 
as an maximization algorithm. Therefore, we take the negative value of the 
objective as the actual output to the implementation. 
%
%
\section{Numerical experiments}\label{se:exp}
\par\noindent
In this section, Algorithms~\ref{alg:descent},~\ref{alg:lv} and~\ref{alg:bo} 
are applied to two-dimensional AM optimization problems in conjunction with the
material model and the finite element method described in 
Section~\ref{se:goveqn}. All finite element simulations employ standard 4-node
rectangular elements with 
3$\times$3 Gaussian quadrature, 
see~\cite{wang2021} for full details. In particular, two numerical examples 
are simulated and the process parameters~$\ybm$ being optimized 
comprise printing speed, layer thickness and convection coefficient. 
\subsection{Two-dimensional wall after heat dissipation}\label{se:opt1}
\par
The first example involves a two-dimensional rectangular wall with width 
of~$20 mm$ and height of~$10 mm$ in its reference configuration. It is 
discretized by~$n_x=40$ elements along the horizontal direction 
and~$n_y=30$ elements in the vertical direction. The bottom side is subjected 
to Dirichlet boundary conditions with zero displacement and constant 
temperature at the ambient value of~$315 K$. The other three sides 
are traction-free and under convection heat transfer conditions 
with~$\theta_{\infty}=315 K$. The initial temperature of the deposited material 
is set to~$500 K$. To help enforce the boundary conditions, an initial material 
layer is placed at the base of the body with its bottom nodes subjected 
to the Dirichlet boundary conditions described above. Each time step
corresponds to the time it takes to print one full element. The time interval 
between the completion of one layer and the initiation of the next one is set 
to one-half the time needed to print a full layer. The simulation 
terminates~$240 s$ after printing is completed to allow the wall 
to cool down close to the ambient temperature so that the shape error 
of different process parameters can be measured in a comparable manner. 
\par
\begin{figure}
 \centering
  \includegraphics[width=0.85\textwidth]{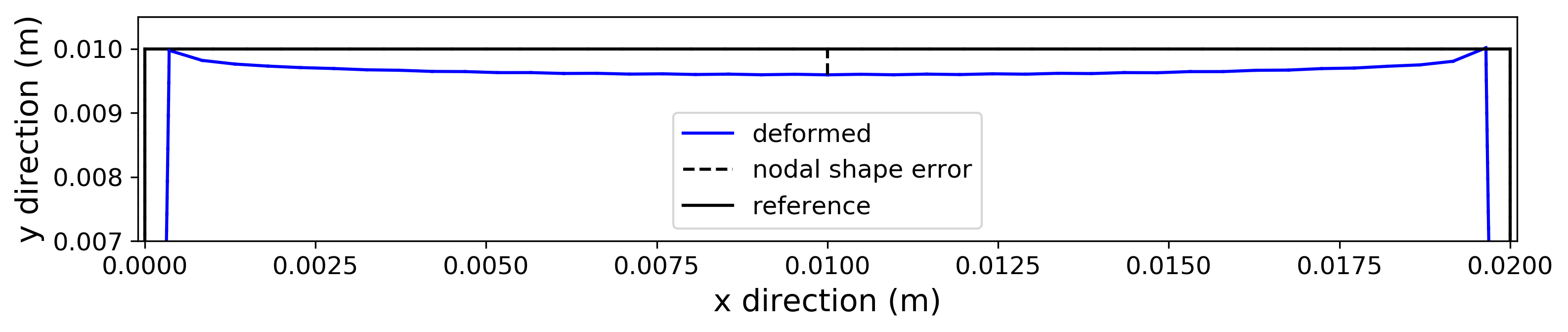}
 \caption{Measurement of shape error for two-dimensional wall 
with deformation magnified five times for visual clarity (the reference
edges in black have zero shape error, while the deformed ones in blue deviate
from the designed height due to the thermomechanical response of the material)}
\label{fig:wall-shape-error}
\end{figure}

For this problem, the shape error is measured only at the top edges 
of the wall, which defines $\Gamma_{S}^{t}$ in~\eqref{eqn:obj-funcl}. 
The local shape 
error~$d\left(\Xbm,\ubm\right)-\bar d(\Xbm)$ in~\eqref{eqn:obj-funcl}
is defined as the difference between the~$y$-coordinate of a 
 point in the deformed configuration and its designed height
at the same~$x$-coordinate of the reference configuration, 
as in Figure~\ref{fig:wall-shape-error}. 
The approximation~\eqref{eqn:obj-funcl-disc} of the surface integral, \textit{i.e.}, the objective, is effected using 3-point Gaussian quadrature
per element edge. Also, the critical length~$L_c$ is chosen to be 
the height of the wall. 
The projected gradient-descent Algorithm~\ref{alg:descent} is applied to the 
optimization of convection coefficient~$h$, which can be controlled through 
mechanisms such as the cooling fan speed in the printing chamber. 
The algorithmic parameters are set to~$\alpha_0 = 1$, $\rho = 0.5$, 
and~$\eta = 0.1$. The time-step size for printing a full element is set 
to~$0.006 s$, equivalent to a printing speed of approximately~$83 mm/s$, 
within the range of FDM. The number of layers is fixed at~$40$. 
The bound constraints on~$h$ are set according to $h\in[30,55]W/m^2K$. 
Clearly, the optimization space here is one-dimensional and the gradient 
$\frac{df_h}{d\ybm}$ can be computed as described in 
Section~\ref{se:grad} and the appendix. The initial value~$h_0$ is set 
to~$40 W/m^2K$. 
\par
\begin{figure}
\centering
\includegraphics[width=0.8\textwidth]{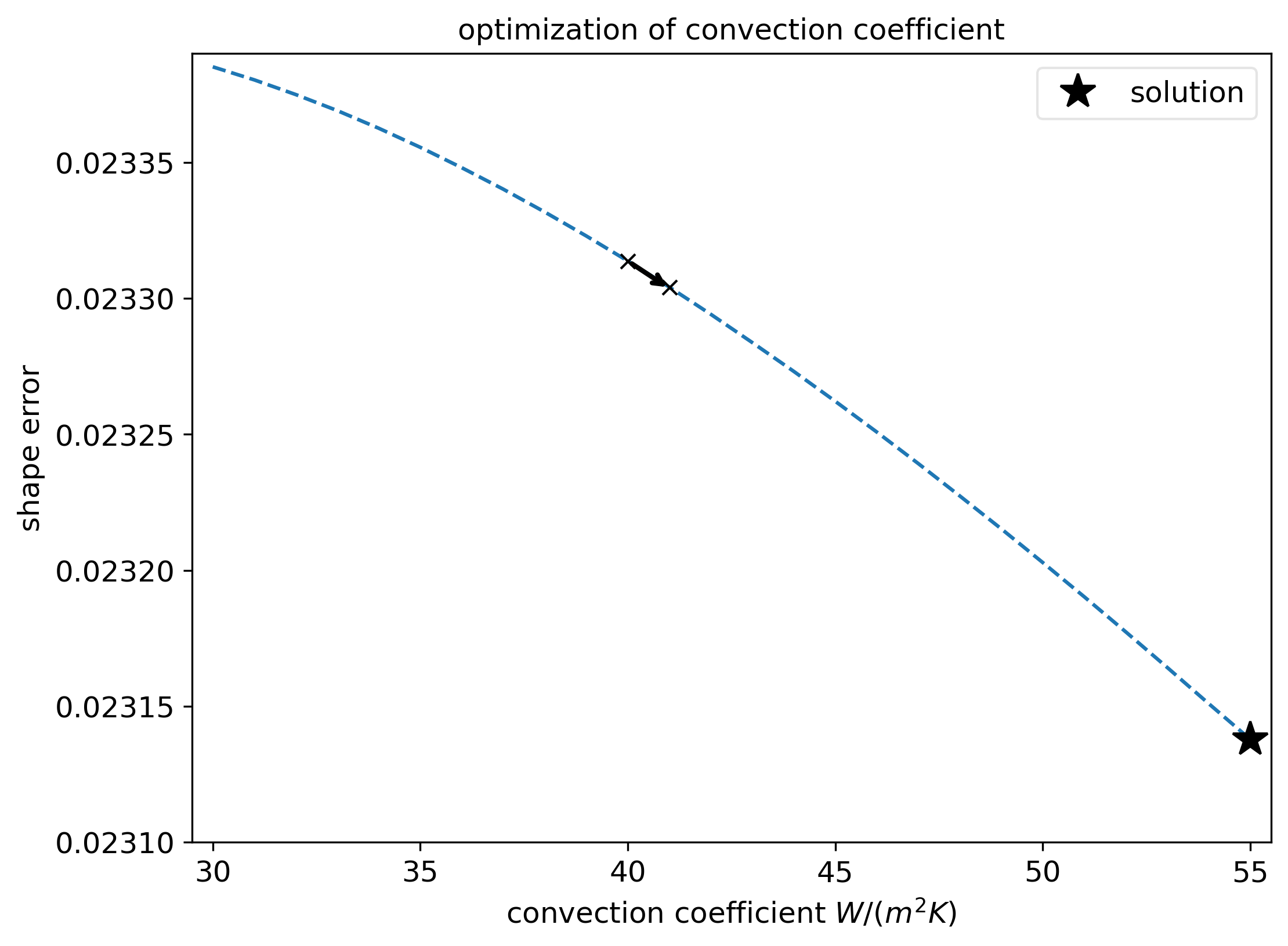}
\caption{Two-dimensional wall: Optimization of convection coefficient
}
\label{fig:1d-opt-convec}
\end{figure}

A plot depicting the results of the optimization is shown 
in Figure~\ref{fig:1d-opt-convec}. Taking a closer look at the initial 
step where $h = 40 W/m^2K$, the objective function is $f_h(40) = 0.0233137$,
while the gradient value is~$\frac{d f_h}{dh}(40)=-9.4\times 10^{-6}$. 
Upon conducting another full simulation for which $f_h(41) = 0.0233041$, 
the step-size $\alpha_k = 1$ satisfies the 
line search condition. Hence, the actual step length is 1 and the next 
convection coefficient in the iteration is~$h=41 W/m^2K$. In a similar manner,
the algorithm eventually reaches the optimal value~$h=55 W/m^2K$ at the 
boundary of the feasible space. The dotted line 
in Figure~\ref{fig:1d-opt-convec} illustrates the monotonically  
decreasing dependence of~$f_h$ on~$h$, which is physically plausible. 
\par
Next, the gradient-free optimization Algorithm~\ref{alg:lv} and~\ref{alg:bo}
are applied to two important process parameters. These are printing speed, 
represented by time-step size~$\delta t$ needed to print one element, 
and layer thickness~$\Delta y$. Clearly, for a given specimen width 
and number of elements in the horizontal direction, the printing speed is
inversely proportional to~$\delta t$. Both process parameters present
challenges when using gradient-based optimization algorithms. 
Layer thickness is a discrete parameter as it is dictated 
by object height and the number of layers. Hence, the objective 
is not differentiable with respect to it. Likewise, printing speed 
is tied to the time-step size~$\Delta t_n$ of the finite element formulation 
in~\eqref{eqn:time-discretization}, and therefore requires additional effort 
for sensitivity calculation. Consequently, gradient-free optimization methods 
present an attractive alternative for these two design parameters. 
\par
The bound constraints for the process parameters are set here to 
\begin{equation}\label{eq:bounds1}
 \centering
  \begin{aligned}
  0.20mm\, \leqslant&\, \Delta y \,\leqslant\, 0.33 mm, \\
  0.005 s\, \leqslant&\, \delta t \,\leqslant\, 0.01 s. \\
  \end{aligned}
\end{equation}
These reflect a realistic range for each parameter. Also, the convection 
coefficient~$h$ is fixed at~$40 W/m^2K$. 
\par
\begin{figure}
  \centering
  \includegraphics[width=0.58\textwidth]{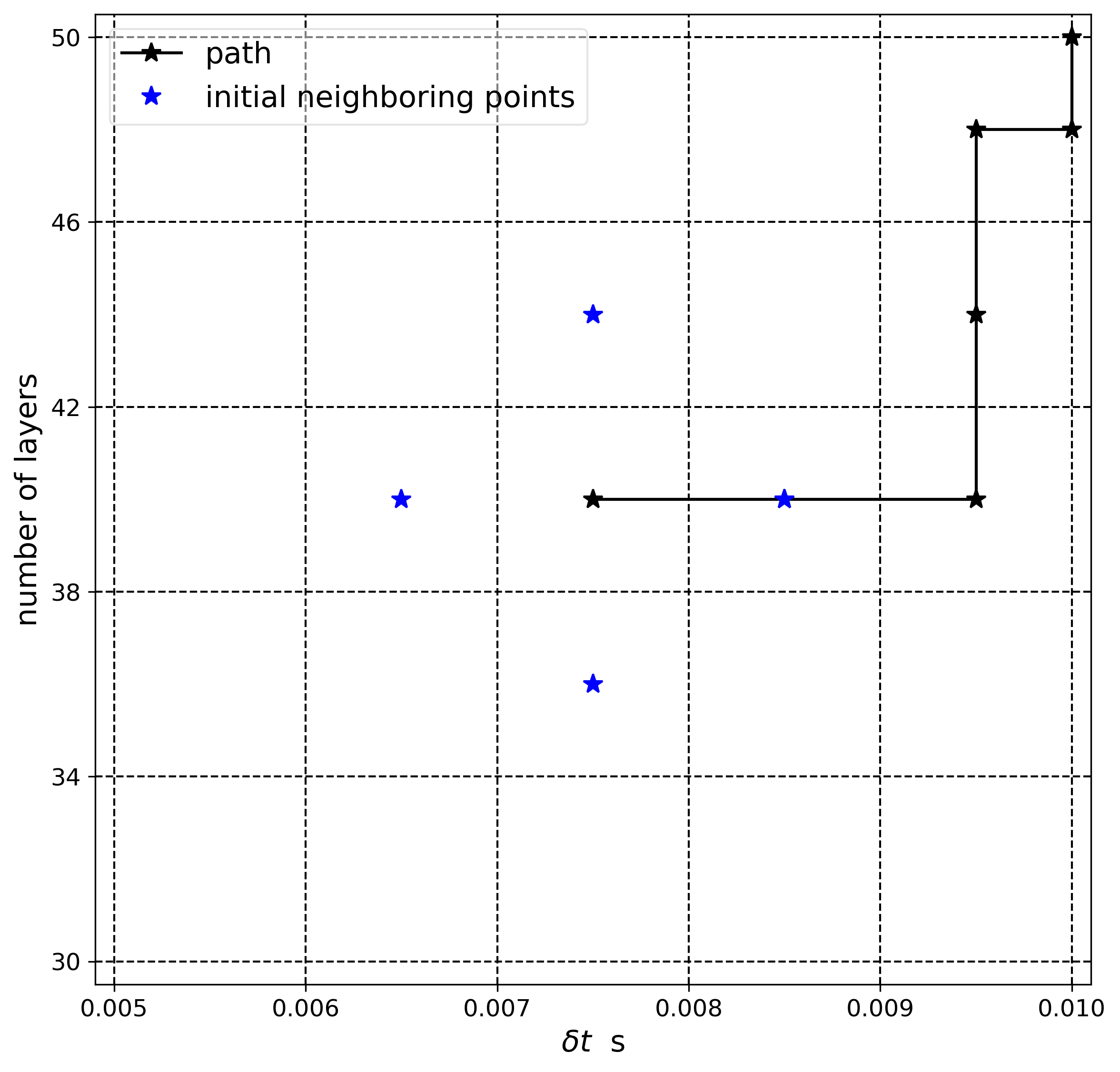}
 \caption{Two-dimensional wall optimization path in the optimization space with method of local variations}
 \label{fig:opt-2d-wall-lv}
\end{figure}

For the method of local variations, the optimization space is a well-defined 
two-dimensional rectangle, as shown in Figure~\ref{fig:opt-2d-wall-lv}. 
The initial optimization step-size $\boldsymbol{\tau}_0$ is set to $4$ for number of
layers and $0.001 s$ for the printing time step. The starting point 
in the parameter space is set to $\Delta y = 0.25 mm$, which corresponds 
to 40 layers, and $\delta t= 0.0075s$, which translates to a printing speed 
of~$v = 76.2 mm/s$. The final step size $\boldsymbol{\tau}_{min}$ is set to be $1$ for the number of layers and $2\times 10^{-4}s$ for the time step.
\par
The optimization path generated by the algorithm follows the rule 
of slower printing speed and larger number of layers and eventually
stops at the boundary governed by the constraints as shown 
in Figure~\ref{fig:opt-2d-wall-lv}. The optimization result translates 
to printing more slowly and slicing the wall more finely in order to reduce 
the shape error, both of which are intuitively plausible. 
Of course, on practical grounds, the trade-off is that the overall printing of
the part takes longer. 
\par
\begin{figure}
 \centering
  \includegraphics[width=0.9\textwidth]{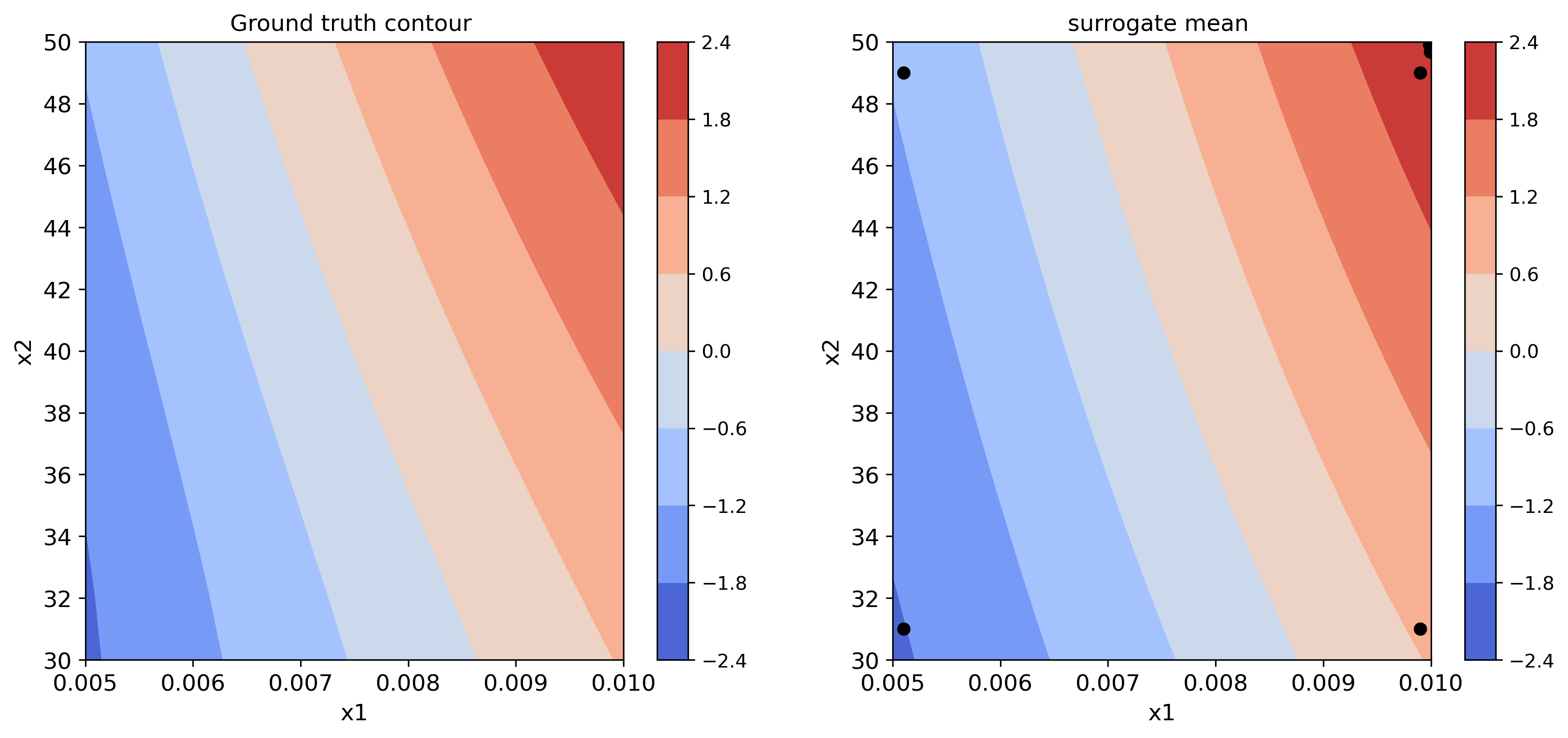}
 \caption{Two-dimensional wall Bayesian optimization, ground truth and surrogate objective.}
\label{fig:wall-bo-obj}
\end{figure}

Note that in this example the purely geometric error, defined 
as the shape error under no deformation, is equal to zero due to the
(nominally) rectangular shape of the printed part. Therefore,  
the thermomechanical response during the building process is responsible 
for the entirety of the shape error. Algorithm~\ref{alg:lv} successfully 
converges to the global minimum for the objective function within the bounds
in~\eqref{eq:bounds1}. The left plot in Figure~\ref{fig:wall-bo-obj} 
is the filled contour plot of the objective function with respect to the 
two process parameters, constructed through exhaustive sampling. 
Evident in the figure is the smoothness of the objective 
function which leads to successful convergence of the optimization procedure. 
Further, the objective function is monotonic in both process parameters within 
the bounds, making it easy to find neighboring points with a smaller objective 
function value.  
\par
\begin{figure}
 \centering
  \includegraphics[width=0.9\textwidth]{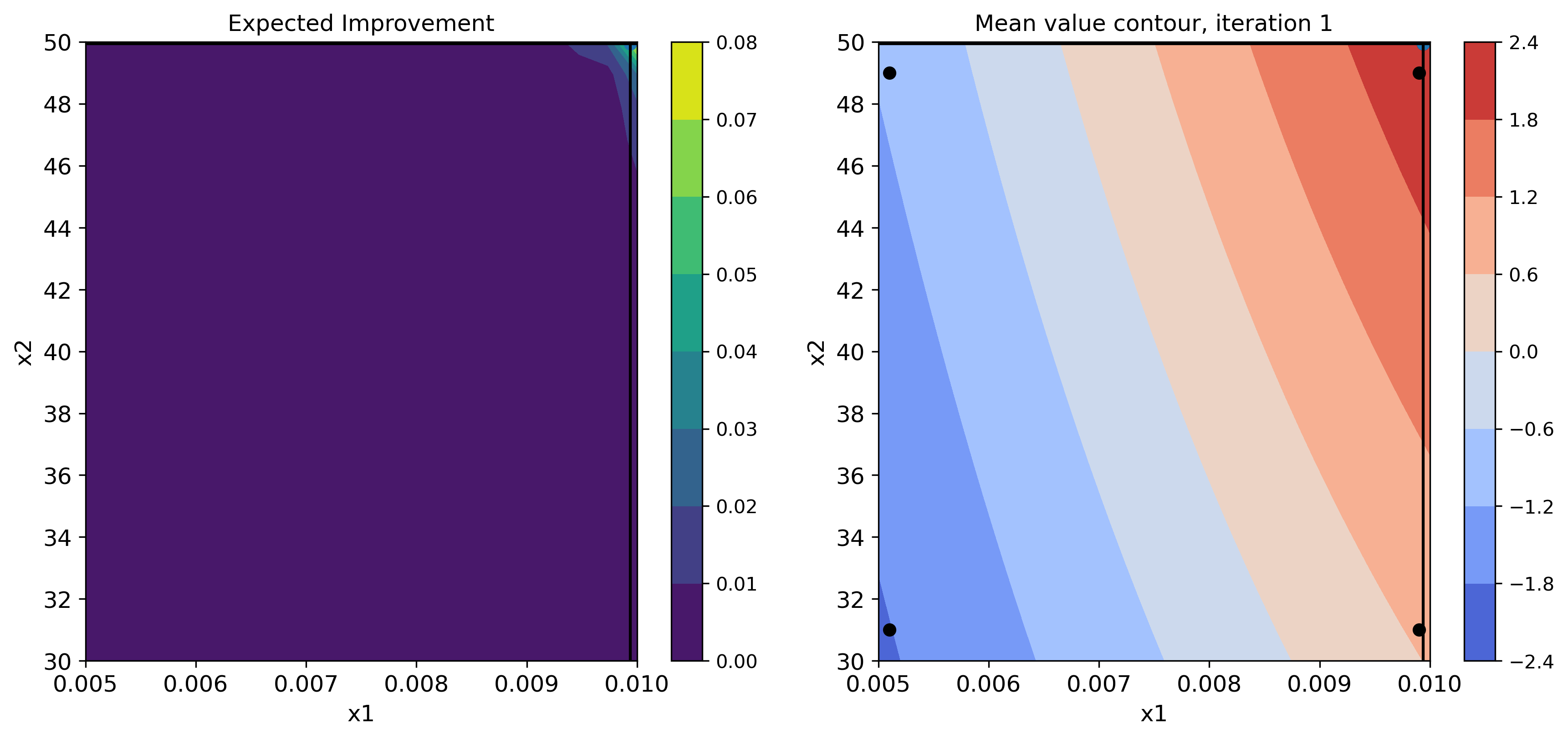}
 \caption{Two-dimensional wall Bayesian optimization, expected improvement and mean value contour, iteration 1. 
         The axes $x_1$ and $x_2$ are the number of layers ($10/\Delta y$) and temporal step size $\delta t$,
respectively, as in~\eqref{eq:bounds1}.}
\label{fig:wall-bo-it-1}
\end{figure}
\begin{figure}
 \centering
  \includegraphics[width=0.9\textwidth]{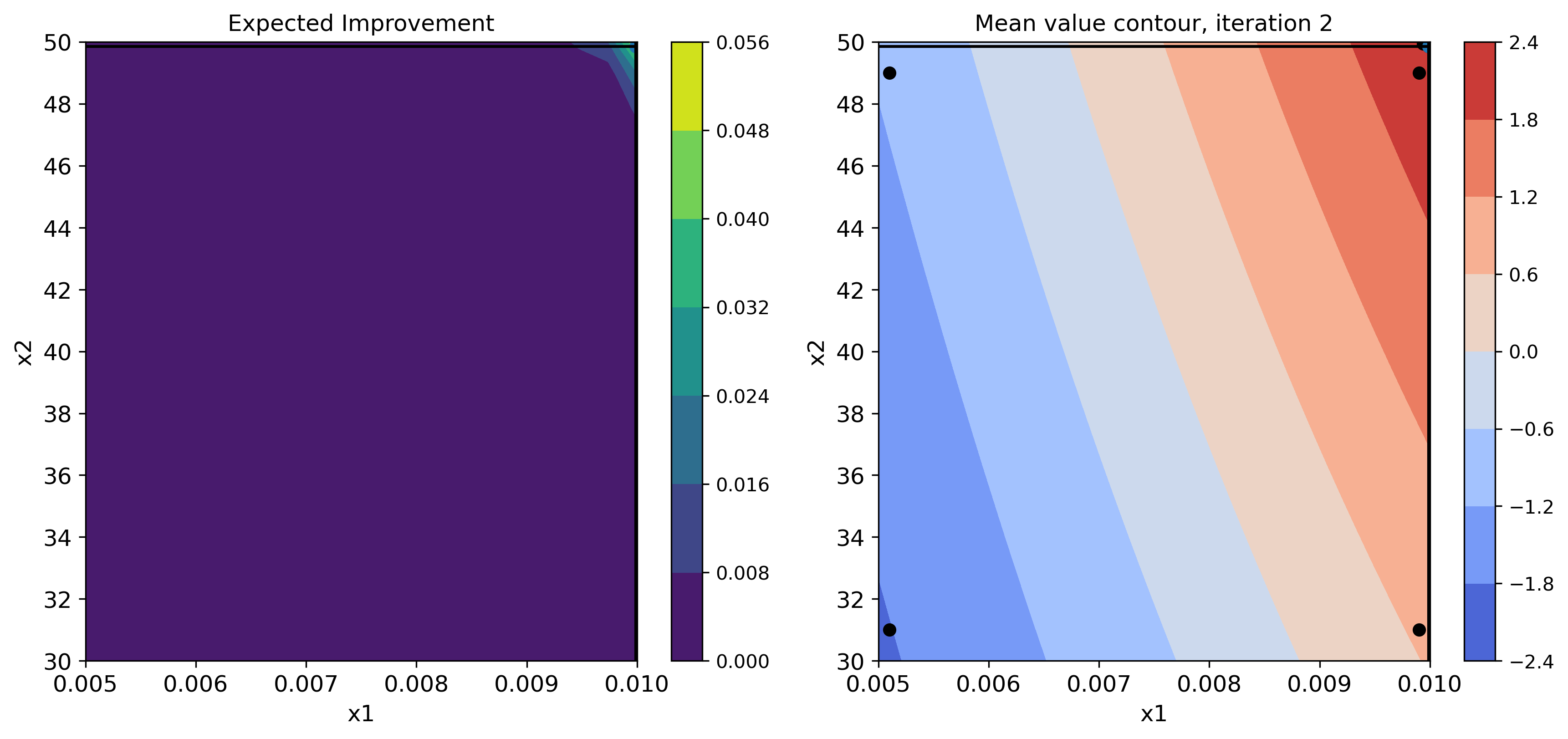}
 \caption{Two-dimensional wall Bayesian optimization, expected improvement and mean value contour, iteration 2. The black dots represent sampled points.}
\label{fig:wall-bo-it-2}
\end{figure}

The initial sampling is important in the convergence of the Bayesian optimization. 
Here, four corner points are selected in the bounded optimization space. 
The mean value and expected improvement acquisition function value contour for each iteration are shown in 
Figure~\ref{fig:wall-bo-it-1} and~\ref{fig:wall-bo-it-2}. 
The ground truth objective 
and the prediction by the surrogate model at termination are displayed in 
Figure~\ref{fig:wall-bo-obj}, where the black dots in the surrogate contour are 
the sampled points. Notice that the objective value plotted is the scaled value 
of the $-f_h$ since Algorithm~\ref{alg:bo} solves a equivalent maximization problem, 
as explained in Section~\ref{se:bayesian}. 
Thus, the maximum value shown indeed corresponds to the minimal shape error. In 
this case, the four initial points suffice in generating a decent surrogate 
model and the algorithm terminates in 2 iterations as the predicted optimal 
variables produced by two consecutive iterations are close enough, at the top-right corner in Figure~\ref{fig:wall-bo-obj}. 
The optimal design variables are found to be~$50$ layers and~$0.01 s$, which 
translates to layer thickness of~$0.2 mm$ and a printing speed of~$0.057m/s$. 
%
\subsection{Two-dimensional wall with hole after heat dissipation}\label{se:opt2}
\par
\begin{figure}
  \centering
  \includegraphics[width=0.7\textwidth]{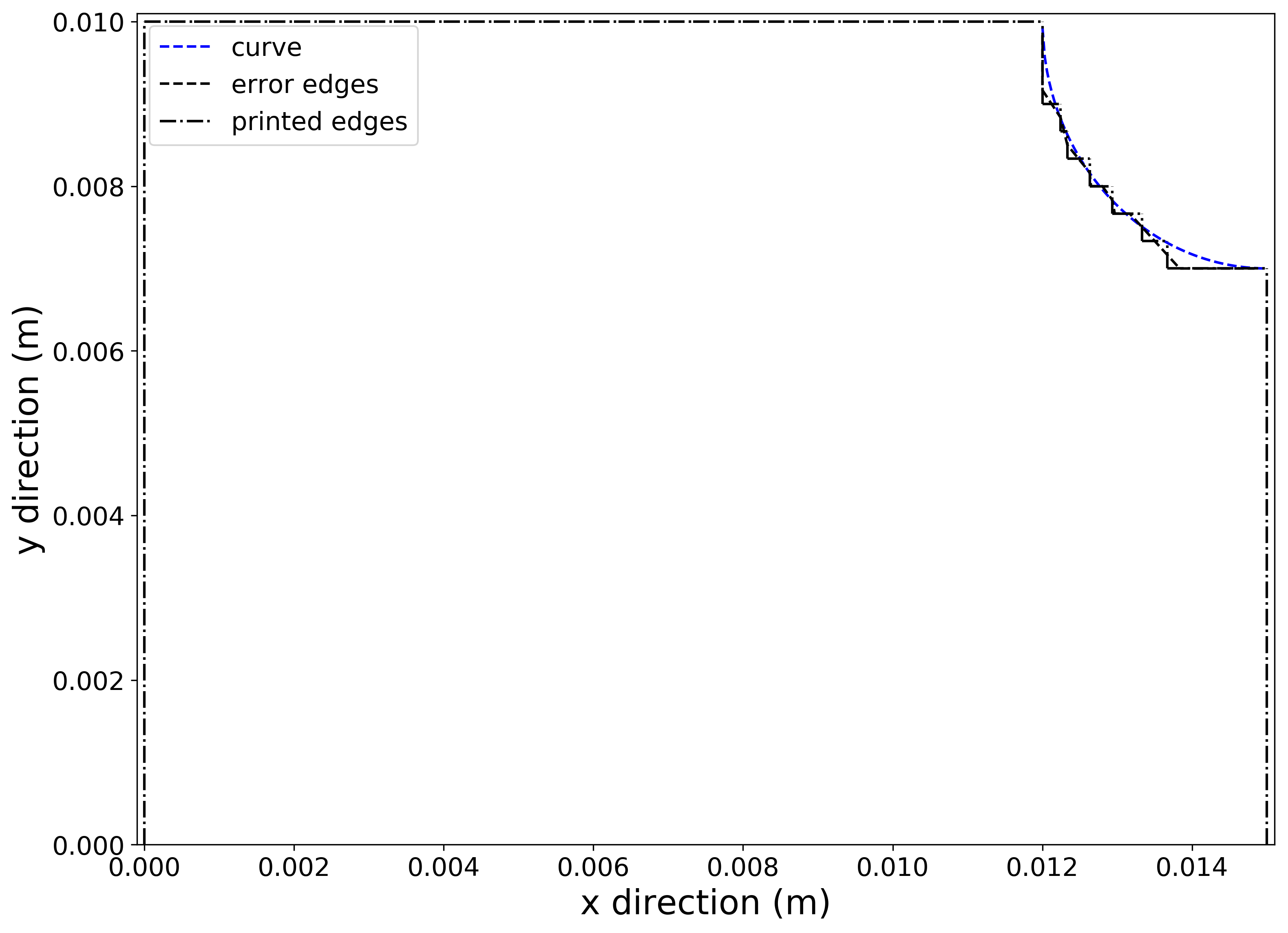}
 \caption{Boundary and error edges for two-dimensional wall with hole}
 \label{fig:error-edge}
\end{figure}

This example concerns the simulation of printing a two-dimensional wall 
with a circular hole. The width of the wall is~$15 mm$ and the height 
is~$10 mm$. A quarter-circular hole of radius~$r=3 mm$ is situated 
at the top-right corner of the wall, as in Figure~\ref{fig:error-edge}.  
Due to the curved shape, hanging nodes are used in 
the vicinity of the circular boundary mesh. 
\par
\begin{figure}
  \centering
  \includegraphics[width=0.7\textwidth]{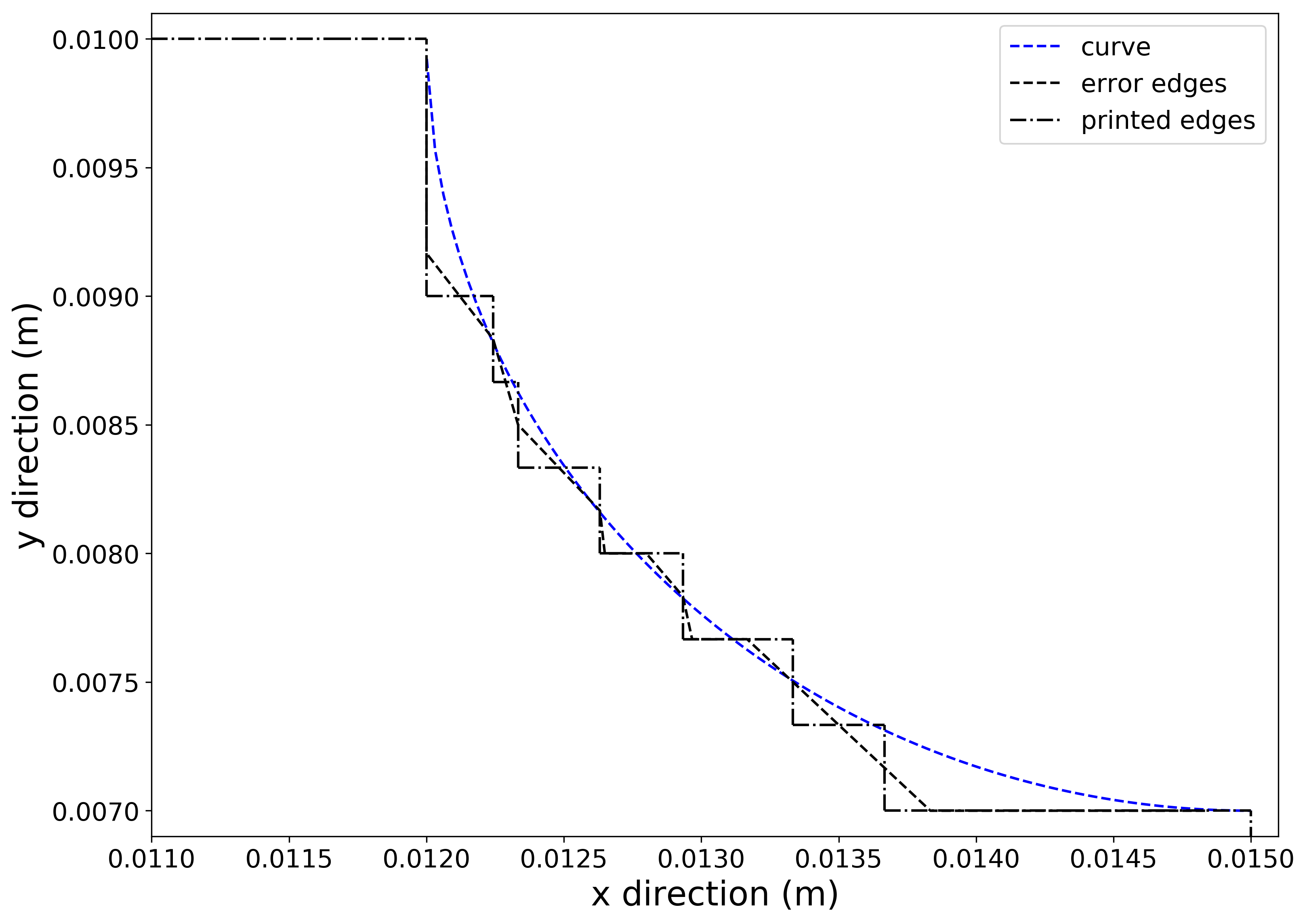}
 \caption{Error edges for two-dimensional wall with hole zoomed in}
 \label{fig:error-edge-zoom}
\end{figure}

The shape error is measured along the circular edge. Specifically, 
$\Gamma_{S}^{t}$ in~\eqref{eqn:obj-funcl} is the {\em step-edge}  
approximation of the quarter-circular hole in the reference configuration,
which is defined as the edges formed by connecting the mid-points 
of the actual printed-edge approximation of the hole, 
as in Figure~\ref{fig:error-edge-zoom}. This definition ensures that 
in the reference configuration the length of~$\Gamma_{S}^{t}$, 
under mesh refinement, becomes closer to the arc length of the circle~$\pi r/2$, in contrast to the length 
of the actual mesh boundary, which is always equal to~$2r$. 
The local shape error~${d\left(\Xbm,\ubm\right)-\bar d(\Xbm})$ 
in~\eqref{eqn:obj-funcl-disc} is defined as the deviation in radius 
of the step-edge approximation from the circle of the designed shape. 
The critical length~$L_c$ here is chosen to be the radius~$r$ of the 
same circle.
\par
\begin{figure}
\centering
\includegraphics[width=0.8\textwidth]{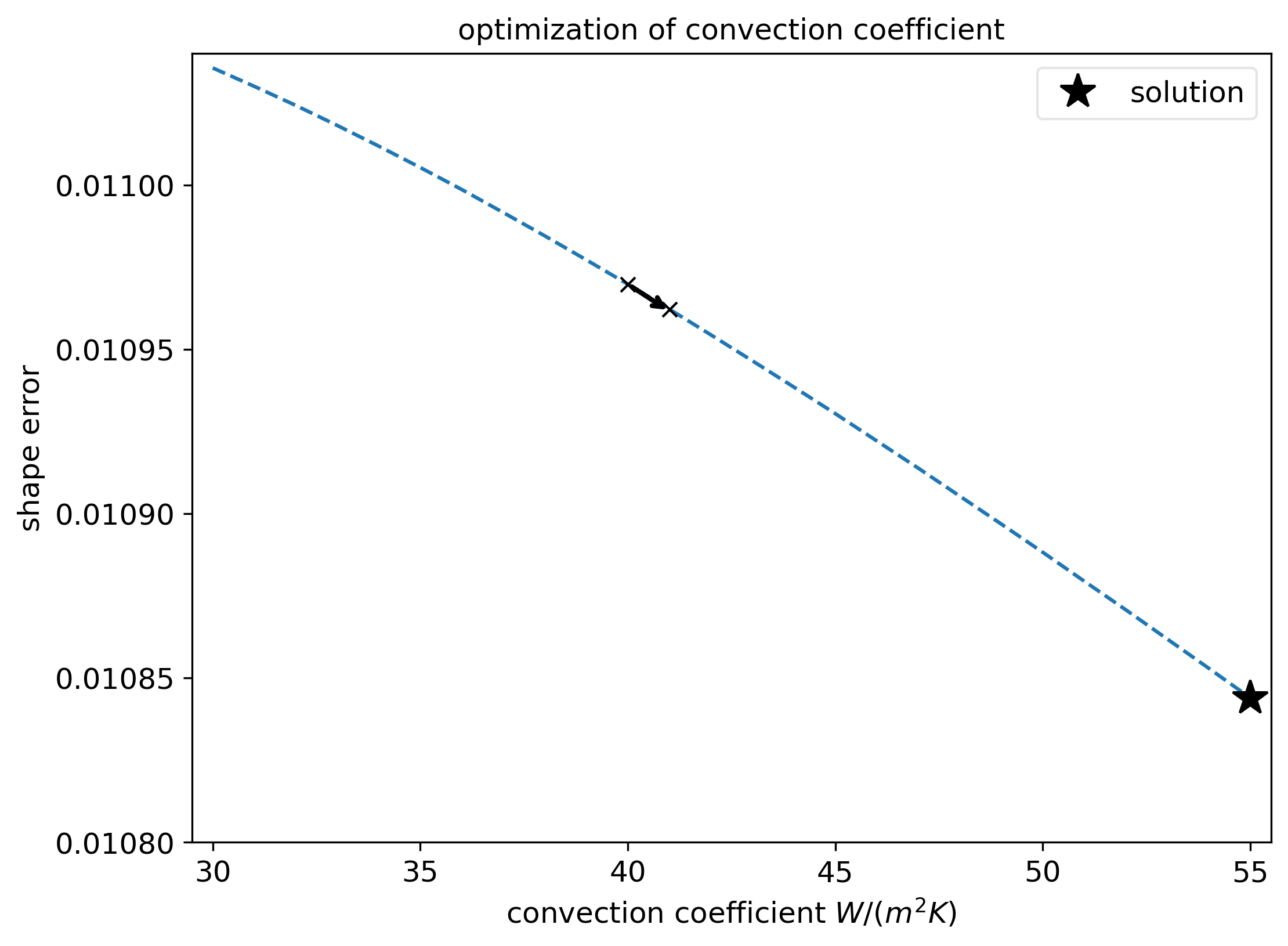}
\caption{Two-dimensional wall with hole: Optimization of convection coefficient}
\label{fig:1d-cir-opt-convec}
\end{figure}

Algorithm~\ref{alg:descent} is first applied 
to determine the optimal convection coefficient $h\in[30,55]W/m^2K$, 
where~$ \alpha_0 = 1$, $\rho = 0.5$ and~$\eta = 0.1$. 
The time-step size for printing a full element is set 
to~$0.006 s$, equivalent to a printing speed of~$62.5 mm/s$, while 
the number of layers is fixed at~$30$. Again, the 
sensitivity~$\frac{df_h}{d\ybm}$ is computed based on the method described 
in Section~\ref{se:grad} and the appendix. However, additional care needs 
to be taken due to the presence of hanging nodes around the quarter-circle. 
As described towards the end of Section~\ref{se:fem}, a ratio is maintained for the displacement and temperature of a hanging node that resides 
on an edge with two other nodes at the end of the edge. 
Naturally, the sensitivity of the displacement and temperature of the hanging node is required to maintain the same ratio.
The initial convection 
coefficient~$h_0$ is set to~$40 W/m^2K$. 
The optimization path is shown in Figure~\ref{fig:1d-cir-opt-convec}. 
%
\par
The two gradient-free optimization algorithms are applied to the optimization 
of the printing speed and layer thickness for~$h=40 W/m^2K$. 
The bounds of the process parameters are set to 
\begin{equation}
 \centering
  \begin{aligned}
  0.20 mm\, \leqslant\,& \Delta y\, \leqslant\, 0.33 mm\ , \\
  0.004 s\, \leqslant\,& \delta t\, \leqslant\, 0.009 s\ . \\
  \end{aligned}
\end{equation}
The remaining algorithmic parameter setup is as in the previous example. 
However, the presence of the circular hole causes the shape error 
to have a more complex relation to the number of layers. 
Indeed, it is evident that the shape error here is nonzero in the reference configuration with no deformation.
This is because the shape error is measured against the ideal quarter-circle 
on rectangular boundary elements. The referential shape error is intrinsic to the layered
approximation of a circular curve. 
The shape error is also affected by the slicing algorithm that determines where each layer should 
stop at the boundary of the circle. 
The discrete nature of the number of layers, slicing algorithm and rectangular approximation of a circle lead to a 
non-monotonic relationship between the shape error and the number of layers. We refer to the shape error in 
reference configuration with no deformation as the geometric error.
The geometric error does not change the general trend of how the process parameters affect the shape error, as we show later in the section,
because the thermomechanical shape error caused by the AM building process itself is significantly larger.
More discussion on this topic can be found in~\cite{wang2021,choi2016influence,paul2014effect}.
\par
\begin{figure}
  \centering
  \includegraphics[width=0.58\textwidth]{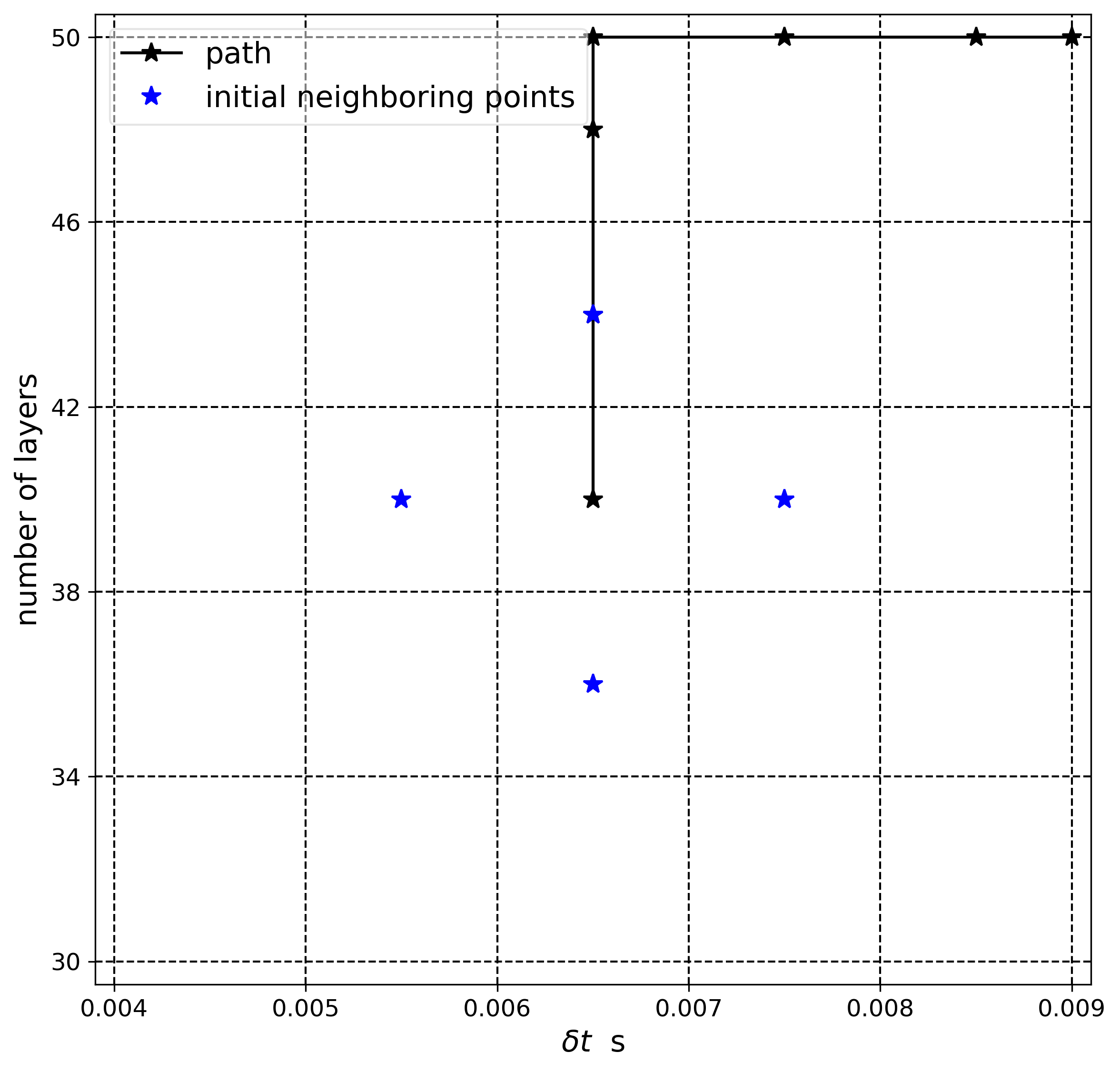}
 \caption{Two-dimensional wall with hole: Optimization path with method of local variations}
 \label{fig:opt-2d-quacir-lv}
\end{figure}
\begin{figure}
 \centering
  \includegraphics[width=0.9\textwidth]{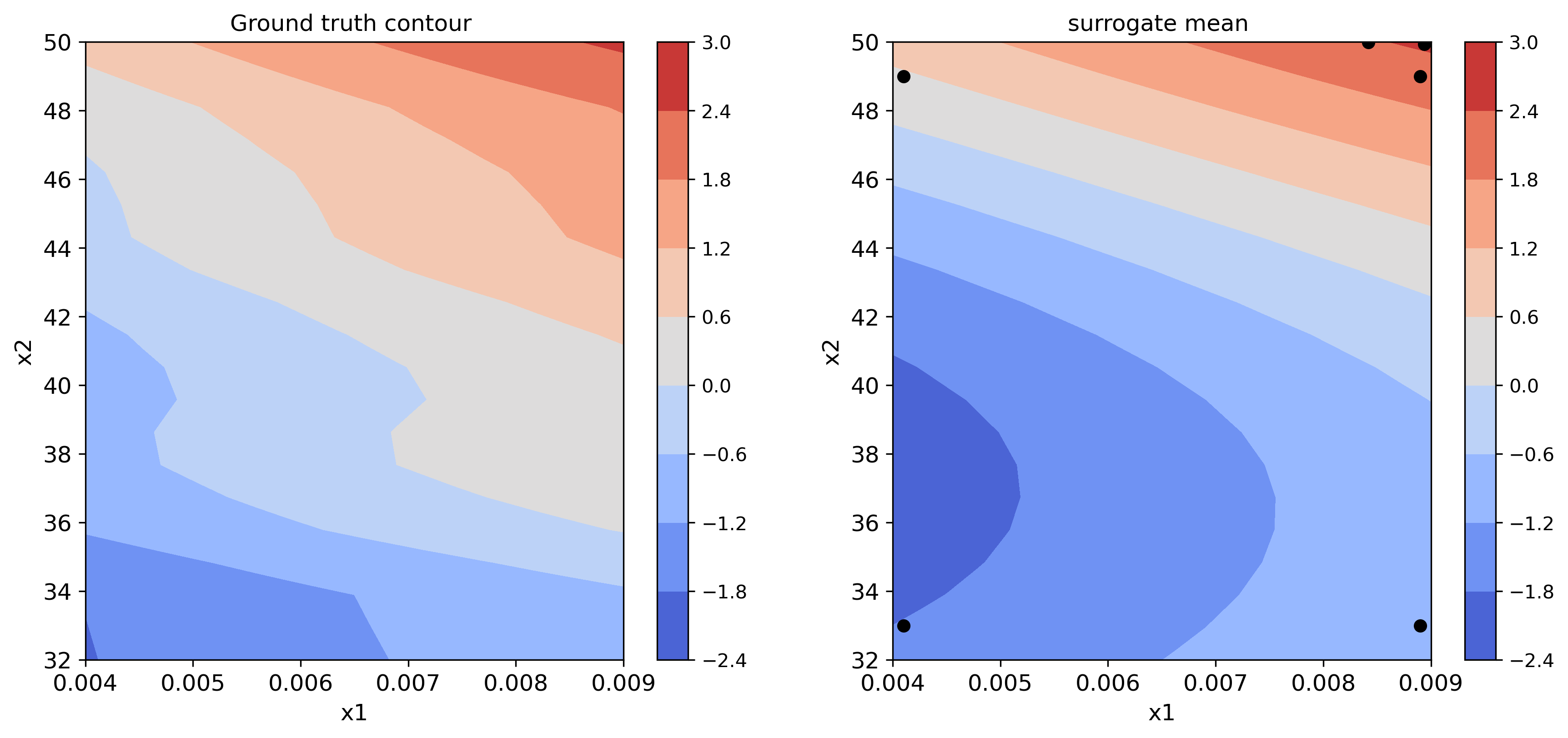}
 \caption{Two-dimensional wall with hole: Bayesian optimization: Ground truth 
and surrogate objective}
\label{fig:wall-cir-bo-obj}
\end{figure}

The optimization path generated by the method of local variations 
is shown in Figure~\ref{fig:opt-2d-quacir-lv}, where the initial
step size is~$0.001 s$ for the print velocity and~$4$ for the number of layers. 
Starting from the point $\left(0.0065 s,40\right)$ in 
printing speed/number of layer-space, the optimization method traces a path 
to $\left(0.0065 s,50\right)$ followed by consecutively smaller speed until
the boundary is hit at $\left(0.009 s,50\right)$. The step size $\taubold_k$ 
keeps decreasing to $\taubold_{min}$ (1 for the number of layers) 
as no neighboring point produces smaller objective function values and 
eventually the iteration stops. Therefore, the optimal variables given by the method of local variations are $\left(0.009 s,50\right)$,
which are optimal as shown in the filled contour plot on the left in Figure~\ref{fig:wall-cir-bo-obj}. 

\par
\begin{figure}
 \centering
  \includegraphics[width=0.9\textwidth]{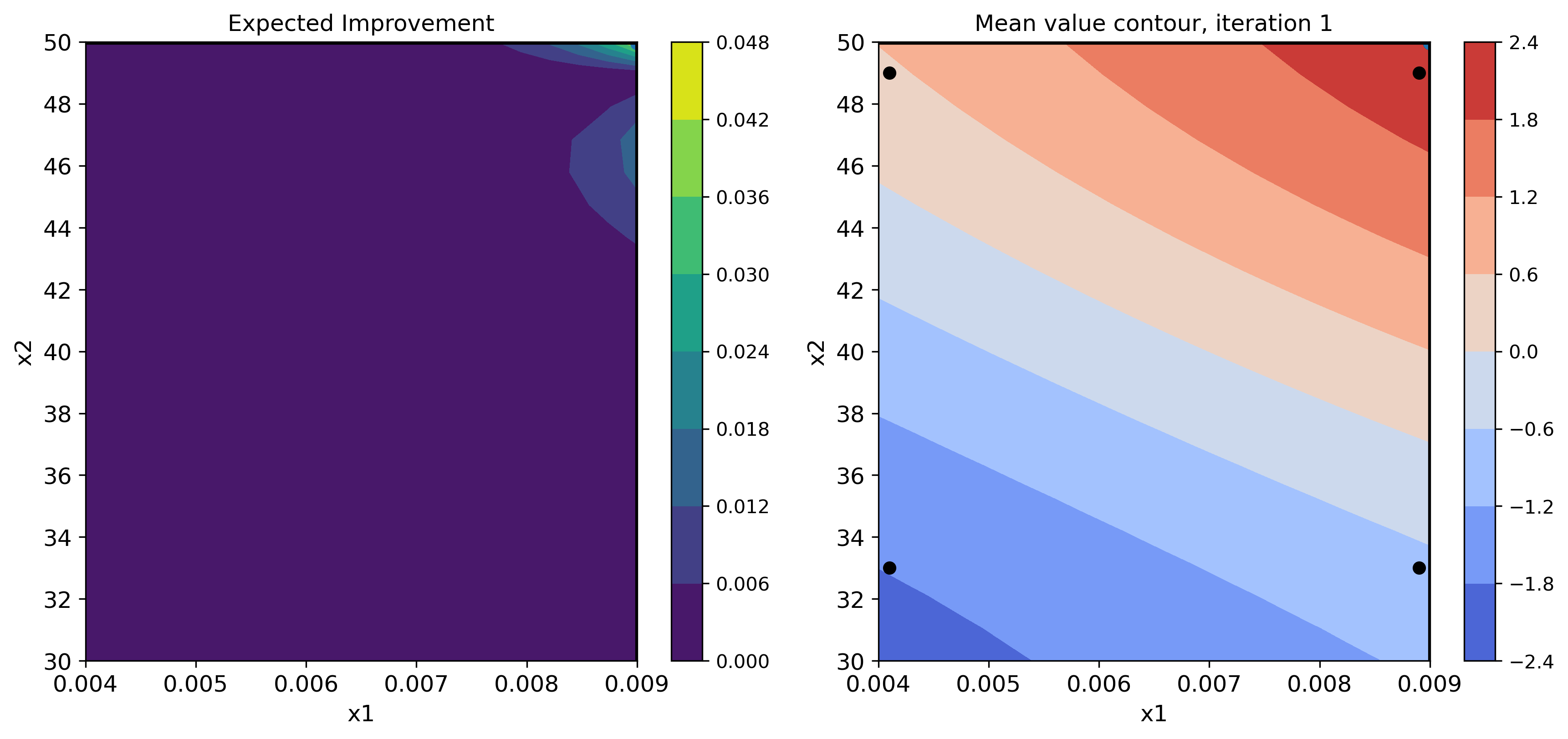}
 \caption{Two-dimensional wall with hole: Bayesian optimization, expected improvement and mean value contour, iteration 1. The black dots are the sampled points.}
\label{fig:wall-cir-bo-it-1}
\end{figure}
\begin{figure}
 \centering
  \includegraphics[width=0.9\textwidth]{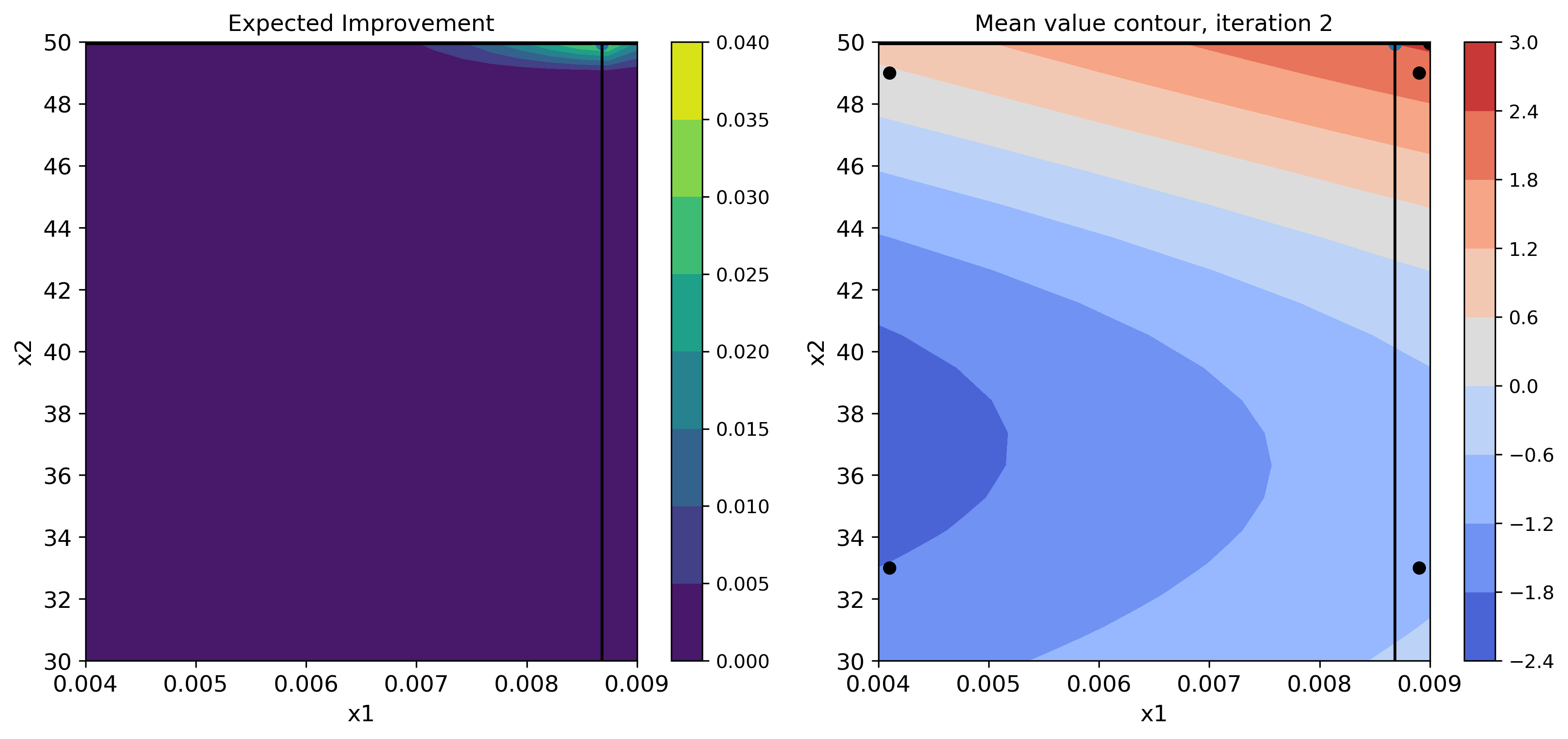}
 \caption{Two-dimensional wall with hole Bayesian optimization, expected improvement and mean value contour, iteration 2. The intersection of horizontal and vertical lines are the new sampled points to be included in the surrogate model.}
\label{fig:wall-cir-bo-it-2}
\end{figure}

Next, the Bayesian optimization algorithm is applied to the same problem, 
starting from four initial sample points. The mean value and expected improvement contour 
for selected iterations are shown in Figures~\ref{fig:wall-cir-bo-it-1} 
and~\ref{fig:wall-cir-bo-it-2}, respectively. 
The ground truth objective and the prediction by the surrogate model 
at the time of termination are displayed in Figure~\ref{fig:wall-cir-bo-obj}, 
where the black dots in the surrogate contour are the sampled points. 
It is emphasized again that the objective shown in the plots are 
the scaled values of negative shape error. 
\par
\begin{figure}
 \centering
  \includegraphics[width=0.9\textwidth]{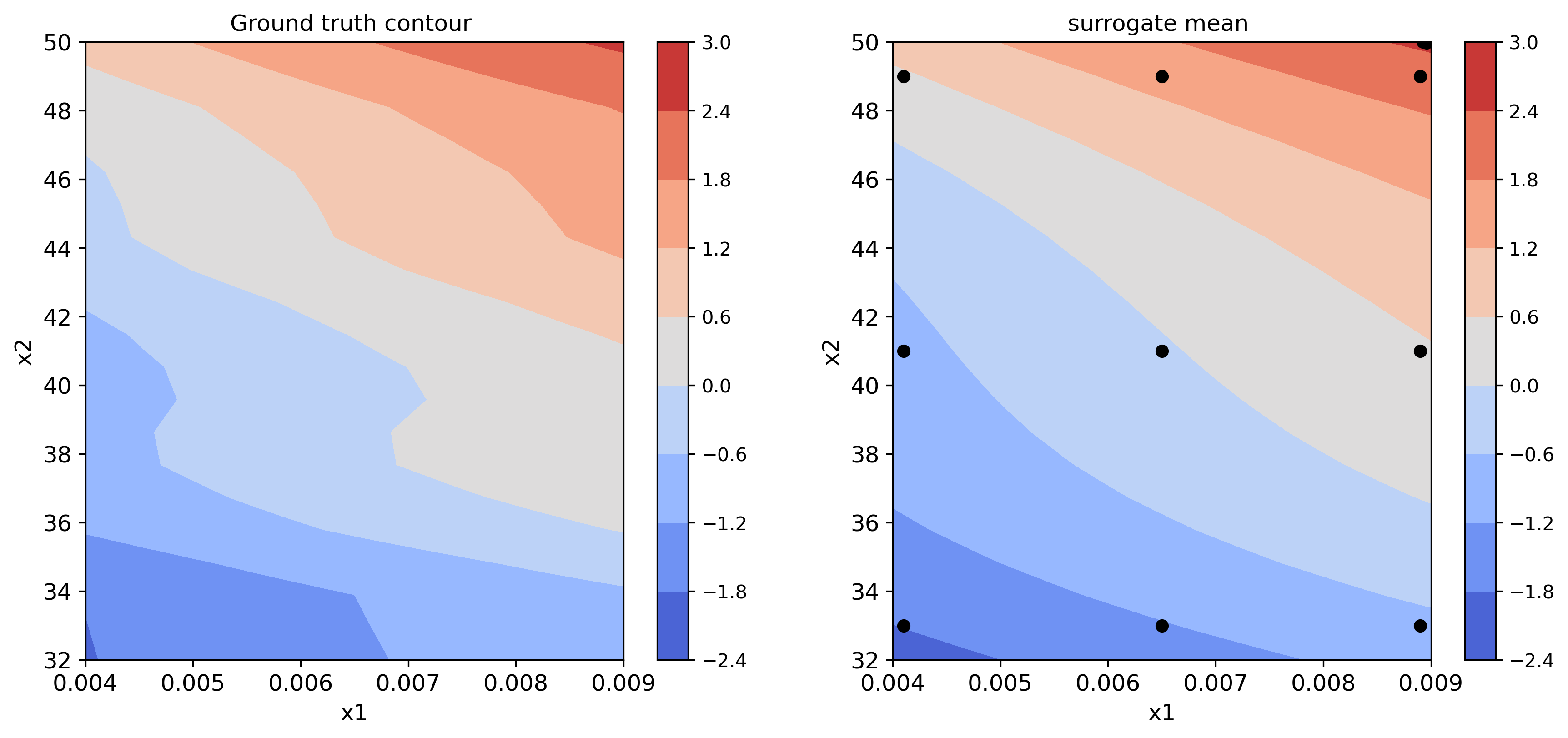}
 \caption{Two-dimensional wall with hole: Bayesian optimization, 9 initial samples, ground truth and surrogate objective}
\label{fig:wall-cir-bo-9-obj}
\end{figure}

While Algorithm~\ref{alg:bo} still converges in two iterations, it is clear 
that the predicted mean value contour is not an accurate depiction of the ground truth, as shown in Figure~\ref{fig:wall-cir-bo-obj}. However, since the algorithm 
only strives to find the optimal design variables, it terminates as the point 
$(0.009 s, 50)$ at the top right corner is found repeatedly. If we have 
nine initial sample points, as shown in Figure~\ref{fig:wall-cir-bo-9-obj}, an 
overall more accurate surrogate model for~$f_h$ of~$\Delta y$ and~$\delta t$ can be 
found. 
\par
The results from both examples predict a better finish quality with higher 
cooling rate (convection coefficient) within a realistic range. They 
both point to reduced shape error with slower printing and more layers. These 
seemingly obvious results would be considerably more complicated given real 
application constraints such as production time and cost associated with 
increased number of layers. In addition, the wall with a circular hole example 
illustrates that both the process parameters studied in this work and the pure 
geometric design of the object contribute to the shape error. The 
gradient-based optimization, while costly, is more accurate and can scale well. 
It can also rely on the wealth of optimization theories to guarantee local 
minimum convergence~\cite{optimization}. However, computing the sensitivities takes 
effort and is in general less flexible to new variables and problems. 
\par
Gradient-free optimization is a promising alternative, particularly Bayesian optimization. 
This approach eliminates the need to compute 
complicated sensitivities, and thus can be easily adapted to different simulation models, including ones different from ours. Accompanied by physical 
experiments for validation, it can be useful in practice. In this work, we rely 
on a fully coupled thermomechanical finite element model to generate sample points. In 
practice, lower fidelity weakly coupled models can work well 
and further reduce the cost of each simulation. Bayesian optimization can be 
applied directly to those more commonly used production models. We point out 
that if the optimization space is large, the method using GP might not scale well. 
%
%
\section{Conclusions}\label{se:conclusion}
\par\noindent
In this paper, design optimization methods for AM process parameters have been 
proposed based on finite element simulation. To address the diverse set of 
variables, both gradient-based and gradient-free optimization methods 
have been proposed and implemented. While the objective is chosen to be 
shape error, the framework can be adapted to other design goals. The 
gradient-based approach relies on fully discretized reduced space formulation 
of the balance equations to obtain sensitivities. For 
gradient-free optimization, a local search algorithm and a Bayesian optimization algorithm using Gaussian process surrogate models and expected improvement acquisition functions are presented. Numerical 
experiments demonstrate that these optimization methods can yield physically 
plausible results efficiently. Line search gradient descent algorithm is applied to 
convection coefficient and gradient-free algorithms are deployed for layer 
thickness and printing speed. These flexible methods can be easily extended to 
other design variables. And in the case of gradient-free optimization,  
other finite element models, particularly mature uncoupled three-dimensional 
ones, can be easily adopted. Surrogate models using NN can also be developed.
%
%
\small
\addcontentsline{toc}{section}{References}
\lsp\bibliography{paper2}

\begin{thebibliography}{10}

\bibitem{mohamed2015}
O.A. Mohamed, S.H. Masood, and J.L. Bhowmik.
\newblock Optimization of fused deposition modeling process parameters: A
  review of current research and future prospects.
\newblock {\em Advances in Manufacturing}, 3:42--53, March 2015.

\bibitem{gibson2014}
I.~Gibson, D.W. Rosen, and B.~Stucker.
\newblock {\em Additive Manufacturing Technologies}.
\newblock Springer, New York, 2014.

\bibitem{chua2010}
C.K. Chua, K.F. Leong, and C.S. Lim.
\newblock {\em Rapid prototyping: principles and applications (with companion
  CD-ROM)}.
\newblock World Scientific Publishing Company, Singapore, 2010.

\bibitem{upcraft2003}
S.~Upcraft and R.~Fletcher.
\newblock The rapid prototyping technologies.
\newblock {\em Assembly Automation}, 23:318--330, December 2003.

\bibitem{mansour2003}
S.~Mansour and R.~Hague.
\newblock Impact of rapid manufacturing on design for manufacture for injection
  moulding.
\newblock {\em Proceedings of the Institution of Mechanical Engineers, Part B:
  Journal of Engineering Manufacture}, 217:453--461, April 2003.

\bibitem{hopkinson2006}
N.~Hopkinson, R.J. Hague, and P.M. Dickens.
\newblock {\em {Rapid manufacturing: An Industrial Revolution for the Digital
  Age}}.
\newblock John Wiley and Sons, 2006.

\bibitem{thavornyutikarn2014}
B.~Thavornyutikarn, N.~Chantarapanich, K.~Sitthiseripratip, G.A. Thouas, and
  Q.~Chen.
\newblock Bone tissue engineering scaffolding: Computer-aided scaffolding
  techniques.
\newblock {\em Progress in Biomaterials}, 3:61--102, December 2014.

\bibitem{ahn2002}
S.H. Ahn, M.~Montero, D.~Odell, S.~Roundy, and P.K. Wright.
\newblock Anisotropic material properties of fused deposition modeling {ABS}.
\newblock {\em Rapid Prototyping Journal}, 8:248--257, October 2002.

\bibitem{francois2017modeling}
M.M. Francois, A.~Sun, W.E. King, N.J. Henson, D.~Tourret, C.A. Bronkhorst, N.N
  Carlson, C.K. Newman, T.S. Haut, J.~Bakosi J, and J.W. Gibbs.
\newblock Modeling of additive manufacturing processes for metals: Challenges
  and opportunities.
\newblock {\em Current Opinion in Solid State and Materials Science},
  21:SAND--2017--6832J, January 2017.

\bibitem{llnl-am-overview}
W.~King, A.T. Anderson, R.M. Ferencz, N.E. Hodge, C.~Kamath, and S.A.
  Khairallah.
\newblock Overview of modelling and simulation of metal powder bed fusion
  process at {Lawrence Livermore National Laboratory}.
\newblock {\em Materials Science and Technology}, 31:957--968, June 2015.

\bibitem{overview2017}
B.~Schoinochoritis, D.~Chantzis, and K.~Salonitis.
\newblock Simulation of metallic powder bed additive manufacturing processes
  with the finite element method: A critical review.
\newblock {\em Proceedings of the Institution of Mechanical Engineers, Part B:
  Journal of Engineering Manufacture}, 231:96--117, January 2017.

\bibitem{hajializadeh2019}
F.~Hajializadeh and A.~Ince.
\newblock Short review on modeling approaches for metal additive manufacturing
  process.
\newblock {\em Material Design and Processing Communications}, 2:e56, March
  2019.

\bibitem{kalita2003}
S.J. Kalita, S.~Bose, H.L. Hosick, and A.~Bandyopadhyay.
\newblock Development of controlled porosity polymer-ceramic composite
  scaffolds via fused deposition modeling.
\newblock {\em Materials Science and Engineering: C}, 23:611--620, October
  2003.

\bibitem{Schoinochoritis2017}
B.~Schoinochoritis, D.~Chantzis, and K.~Salonitis.
\newblock Simulation of metallic powder bed additive manufacturing processes
  with the finite element method: A critical review.
\newblock {\em Proceedings of the Institution of Mechanical Engineers, Part B:
  Journal of Engineering Manufacture}, 231:96--117, January 2017.

\bibitem{Jayanath2018}
S.~Jayanath and A.~Achuthan.
\newblock A computationally efficient finite element framework to simulate
  additive manufacturing processes.
\newblock {\em Journal of Manufacturing Science and Engineering}, 140, April
  2018.

\bibitem{Michaleris2014}
P.~Michaleris.
\newblock Modeling metal deposition in heat transfer analyses of additive
  manufacturing processes.
\newblock {\em Finite Elements in Analysis and Design}, 86:51--60, September
  2014.

\bibitem{YZhang2006}
Y.~Zhang and Y.K. Chou.
\newblock Three-dimensional finite element analysis simulations of the fused
  deposition modelling process.
\newblock {\em Proceedings of the Institution of Mechanical Engineers, Part B:
  Journal of Engineering Manufacture}, 220:1663--1671, October 2006.

\bibitem{patil2013}
N.~Patil, D.~Pal, and B.~Stucker.
\newblock A new finite element solver using numerical eigen modes for fast
  simulation of additive manufacturing processes.
\newblock In {\em 2013 International Solid Freeform Fabrication Symposium}.
  University of Texas at Austin, August 2013.

\bibitem{ding2014}
J.~Ding, P.~Colegrove, J.~Mehnen, S.~Williams, F.~Wang, and P.S. Almeida.
\newblock A computationally efficient finite element model of wire and arc
  additive manufacture.
\newblock {\em The International Journal of Advanced Manufacturing Technology},
  70:227--236, January 2014.

\bibitem{Montevecchi2016}
F.~Montevecchi, G.~Venturini, A.~Scippa, and G.~Campatelli.
\newblock Finite element modelling of wire-arc-additive-manufacturing process.
\newblock {\em Procedia Cirp}, 55:109--114, January 2016.

\bibitem{michaleris2}
E.R. Denlinger, J.~Irwin, and P.~Michaleris.
\newblock Thermomechanical modeling of additive manufacturing large parts.
\newblock {\em Journal of Manufacturing Science and Engineering}, 136, December
  2014.

\bibitem{michaleris3}
E.R. Denlinger, M.~Gouge, J.~Irwin, and P.~Michaleris.
\newblock Thermomechanical model development and in situ experimental
  validation of the laser powder-bed fusion process.
\newblock {\em Additive Manufacturing}, 16:73--80, August 2017.

\bibitem{Roberts2009}
I.A. Roberts, C.J. Wang, R.~Esterlein, M.~Stanford, and D.J. Mynors.
\newblock A three-dimensional finite element analysis of the temperature field
  during laser melting of metal powders in additive layer manufacturing.
\newblock {\em International Journal of Machine Tools and Manufacture},
  49:916--923, October 2009.

\bibitem{Yang2016}
Q.~Yang, P.~Zhang P, L.~Cheng, Z.~Min, M.~Chyu, and A.C. To.
\newblock Finite element modeling and validation of thermomechanical behavior
  of {Ti-6Al-4V} in directed energy deposition additive manufacturing.
\newblock {\em Additive Manufacturing}, 12:169--177, October 2016.

\bibitem{Kolossov2004}
S.~Kolossov, E.~Boillat, R.~Glardon, P.~Fischer, and M.~Locher.
\newblock {3D FE} simulation for temperature evolution in the selective laser
  sintering process.
\newblock {\em International Journal of Machine Tools and Manufacture},
  44:117--123, February 2004.

\bibitem{wang2021}
J.~Wang and P.~Papadopoulos.
\newblock Coupled thermomechanical analysis of fused deposition using the
  finite element method.
\newblock {\em Finite Elements in Analysis and Design}, 197:103607, December
  2021.

\bibitem{masood1996}
S.H. Masood.
\newblock Intelligent rapid prototyping with fused deposition modelling.
\newblock {\em Rapid Prototyping Journal}, 2:24--33, March 1996.

\bibitem{groza2007}
J.R. Groza and J.F. Shackelford.
\newblock {\em Materials Processing Handbook}.
\newblock CRC press, Boca Raton, 2007.

\bibitem{bellehumeur2004}
C.~Bellehumeur, L.~Li, Q.~Sun, and P.~Gu.
\newblock Modeling of bond formation between polymer filaments in the fused
  deposition modeling process.
\newblock {\em Journal of Manufacturing Processes}, 6:170--178, January 2004.

\bibitem{chacon2017}
J.M. Chacón, M.A. Caminero, E.~García-Plaza, and P.J. Núñez.
\newblock Additive manufacturing of {PLA} structures using fused deposition
  modelling: Effect of process parameters on mechanical properties and their
  optimal selection.
\newblock {\em Materials and Design}, 124:143--157, June 2017.

\bibitem{nancharaiah2010}
T.~Nancharaiah, D.R. Raju, and V.R. Raju.
\newblock An experimental investigation on surface quality and dimensional
  accuracy of {FDM} components.
\newblock {\em International Journal on Emerging Technologies}, 1:106--111,
  January 2010.

\bibitem{horvath2007}
D.~Horvath, R.~Noorani, and M.~Mendelson.
\newblock Improvement of surface roughness on abs 400 polymer using design of
  experiments {(DOE)}.
\newblock {\em Materials Science Forum}, 561:2389--2392, 2007.

\bibitem{wang2007}
C.C. Wang, T.W. Lin, and S.S. Hu.
\newblock Optimizing the rapid prototyping process by integrating the {Taguchi}
  method with the gray relational analysis.
\newblock {\em Rapid Prototyping Journal}, 13(5):304--315, October 2007.

\bibitem{thrimurthulu2004}
K.P. Thrimurthulu, P.M. Pandey, and N.V. Reddy.
\newblock Optimum part deposition orientation in fused deposition modeling.
\newblock {\em International Journal of Machine Tools and Manufacture},
  44:585--594, May 2004.

\bibitem{mukherjee2017}
T.~Mukherjee, W.~Zhang, and T.~DebRoy.
\newblock An improved prediction of residual stresses and distortion in
  additive manufacturing.
\newblock {\em Computational Materials Science}, 126:360--372, October 2017.

\bibitem{foroozmehr2016}
A.~Foroozmehr, M.~Badrossamay, E.~Foroozmehr, and S.I. Golabi.
\newblock Finite element simulation of selective laser melting process
  considering optical penetration depth of laser in powder bed.
\newblock {\em Materials and Design}, 89:255--263, January 2016.

\bibitem{west2001}
A.P. West, S.P. Sambu, and D.W. Rosen.
\newblock A process planning method for improving build performance in
  stereolithography.
\newblock {\em Computer-Aided Design}, 33:65--79, January 2001.

\bibitem{Strano2013}
G.~Strano, L.~Hao, R.M. Everson, and K.E. Evans.
\newblock Surface roughness analysis, modelling and prediction in selective
  laser melting.
\newblock {\em Journal of Materials Processing Technology}, 213:589--597, April
  2013.

\bibitem{Onwubolu2014}
G.~Onwubolu and F.~Rayegani.
\newblock Characterization and optimization of mechanical properties of {ABS}
  parts manufactured by the fused deposition modelling process.
\newblock {\em International Journal of Manufacturing Engineering}, pages
  1--13, November 2014.

\bibitem{espin2015}
M.~Domingo-Espin, J.M. Puigoriol-Forcada, A.A. Garcia-Granada, J.~Llumà,
  S.~Borros, and G.~Reyes.
\newblock Mechanical property characterization and simulation of fused
  deposition modeling polycarbonate parts.
\newblock {\em Materials and Design}, 83:670--677, October 2015.

\bibitem{ang2006}
K.C. Ang, K.F. Leong, C.K. Chua, and M.~Chandrasekaran.
\newblock Investigation of the mechanical properties and porosity relationships
  in fused deposition modelling-fabricated porous structures.
\newblock {\em Rapid Prototyping Journal}, 12:100--105, March 2006.

\bibitem{sood2010}
A.K. Sood, R.K. Ohdar, and S.S. Mahapatra.
\newblock Parametric appraisal of mechanical property of fused deposition
  modelling processed parts.
\newblock {\em Materials and Design}, 31:287--295, January 2010.

\bibitem{percoco2012}
G.~Percoco, F.~Lavecchia, and L.M. Galantucci.
\newblock Compressive properties of fdm rapid prototypes treated with a low
  cost chemical finishing.
\newblock {\em Research Journal of Applied Sciences, Engineering and
  Technology}, 4:3838--3842, October 2012.

\bibitem{rayegani2014}
F.~Rayegani and G.C Onwubolu.
\newblock Fused deposition modelling {(FDM)} process parameter prediction and
  optimization using group method for data handling {(GMDH)} and differential
  evolution {(DE)}.
\newblock {\em The International Journal of Advanced Manufacturing Technology},
  73:509--519, July 2014.

\bibitem{masood2010}
S.H. Masood, K.~Mau, and W.Q. Song.
\newblock Tensile properties of processed {FDM} polycarbonate material.
\newblock {\em Materials Science Forum}, 654-656:2556--2559, June 2010.

\bibitem{arivazhagan2011}
A.~Arivazhagan, S.~Masood, and I.~Sbarski.
\newblock Dynamic mechanical analysis of fused deposition modelling processed
  polycarbonate.
\newblock {\em Annual Technical Conference - ANTEC, Conference Proceedings},
  1:950--955, January 2011.

\bibitem{arivazhagan2012}
A.~Arivazhagan and S.H. Masood.
\newblock Dynamic mechanical properties of {ABS} material processed by fused
  deposition modelling.
\newblock {\em International Journal of Engineering Research and Applications},
  2:2009--2014, May 2012.

\bibitem{anitha2001}
R.~Anitha, S.~Arunachalam, and P.~Radhakrishnan.
\newblock Critical parameters influencing the quality of prototypes in fused
  deposition modelling.
\newblock {\em Journal of Materials Processing Technology}, 118:385--388,
  December 2001.

\bibitem{sood2009}
A.K. Sood, R.K. Ohdar, and S.S. Mahapatra.
\newblock Improving dimensional accuracy of fused deposition modelling
  processed part using grey {Taguchi} method.
\newblock {\em Materials and Design}, 30:4243--4252, December 2009.

\bibitem{lee2005}
B.H. Lee, J.~Abdullah, and Z.A. Khan.
\newblock Optimization of rapid prototyping parameters for production of
  flexible abs object.
\newblock {\em Journal of Materials Processing Technology}, 169:54--61, October
  2005.

\bibitem{sood2012}
A.K. Sood, R.K. Ohdar, and S.S. Mahapatra.
\newblock Experimental investigation and empirical modelling of {FDM} process
  for compressive strength improvement.
\newblock {\em Journal of Advanced Research}, 3:81--90, January 2012.

\bibitem{cho2000}
H.S Cho, W.S. Park, B.W. Choi, and M.C. Leu.
\newblock Determining optimal parameters for stereolithography processes via
  genetic algorithm.
\newblock {\em Journal of Manufacturing Systems}, 19:18, 2000.

\bibitem{syrlybayev2021}
D.~Syrlybayev, B.~Zharylkassyn, A.~Seisekulova, A.~Perveen, and D.~Talamona.
\newblock Optimization of the warpage of fused deposition modeling parts using
  finite element method.
\newblock {\em Polymers}, 13(21):3849, January 2021.

\bibitem{Gockel2014}
J.~Gockel, J.~Beuth, and K.~Taminger.
\newblock Integrated control of solidification microstructure and melt pool
  dimensions in electron beam wire feed additive manufacturing of {Ti-6Al-4V}.
\newblock {\em Additive Manufacturing}, 1:119--126, October 2014.

\bibitem{Krol2012}
T.A. Krol, M.F. Zaeh, and C.~Seidel.
\newblock Optimization of supports in metal-based additive manufacturing by
  means of finite element models.
\newblock {\em 23rd Annual International Solid Freeform Fabrication Symposium -
  An Additive Manufacturing Conference, SFF 2012}, pages 707--718, 01 2012.

\bibitem{vasinonta2000}
A.~Vasinonta, J.L. Beuth, and M.L. Griffith.
\newblock A process map for consistent build conditions in the solid freeform
  fabrication of thin-walled structures.
\newblock {\em Journal of Manufacturing Science and Engineering}, 123:615--622,
  August 2000.

\bibitem{yin2022}
J.~Yin, W.~Liu, Y.~Cao, L.~Zhang, J.~Wang, Z.~Li, Z.~Zhao, and P.~Bai.
\newblock {Rapid prediction of the relationship between processing parameters
  and molten pool during selective laser melting of cobalt-chromium alloy
  powder: Simulation and experiment}.
\newblock {\em Journal of Alloys and Compounds}, 892:162200, February 2022.

\bibitem{vastola2016}
G.~Vastola, G.~Zhang, Q.X. Pei, and Y.W. Zhang.
\newblock Controlling of residual stress in additive manufacturing of
  {Ti-6Al-4V} by finite element modeling.
\newblock {\em Additive Manufacturing}, 12:231--239, October 2016.

\bibitem{nickel2001}
A.H. Nickel, D.M. Barnett, and F.B. Prinz.
\newblock Thermal stresses and deposition patterns in layered manufacturing.
\newblock {\em Materials Science and Engineering: A}, 317:59--64, October 2001.

\bibitem{alafaghani2017}
A.~Alafaghani, A.~Qattawi, B.~Alrawia, and A.~Guzmana.
\newblock Experimental optimization of fused deposition modelling processing
  parameters: A design-for-manufacturing approach.
\newblock {\em Procedia Manufacturing}, 10:791--803, December 2017.

\bibitem{zhang2006}
Y.~Zhang and Y.K. Chou.
\newblock {3D} {FEA} simulations of fused deposition modeling process.
\newblock In {\em ASME 2006 International Manufacturing Science and Engineering
  Conference}, pages 1121--1128, January 2006.

\bibitem{morgan2016}
H.~D. Morgan, J.~A.Cherry, S.~Jonnalagadda, D.~Ewing, and J.~Sienz.
\newblock Part orientation optimisation for the additive layer manufacture of
  metal components.
\newblock {\em The International Journal of Advanced Manufacturing Technology},
  86:1679--1687, September 2016.

\bibitem{wang2004}
J.~Wang and N.~Zabaras.
\newblock A {Bayesian} inference approach to the inverse heat conduction
  problem.
\newblock {\em International Journal of Heat and Mass Transfer}, 47:3927--3941,
  August 2004.

\bibitem{mathern2021}
A.~Mathern, O.~S. Steinholtz, A.~Sjöberg, et~al.
\newblock Multi-objective constrained bayesian optimization for structural
  design.
\newblock {\em Structural and Multidisciplinary Optimization}, 63:689--701,
  February 2021.

\bibitem{calandra2016}
R.~Calandra, A.~Seyfarth, J.~Peters, and M.~P. Deisenroth.
\newblock Bayesian optimization for learning gaits under uncertainty.
\newblock {\em Annals of Mathematics and Artificial Intelligence}, 76:5--23,
  February 2016.

\bibitem{bernardo2011}
J.~Bernardo, M.~J. Bayarri, J.~O. Berger, A.~P. Dawid, D.~Heckerman, A.~F.
  Smith, and M.~West.
\newblock Optimization under unknown constraints.
\newblock {\em Bayesian Statistics}, 9:229, October 2011.

\bibitem{zuluaga2013}
M.~Zuluaga, G.~Sergent, A.~Krause, and M.~Püschel.
\newblock Active learning for multi-objective optimization.
\newblock In {\em International Conference on Machine Learning}, pages
  462--470. PMLR, February 2013.

\bibitem{zhang2021}
M.~Zhang, A.~Parnell, D.~Brabazon, and A.~Benavoli.
\newblock Bayesian optimisation for sequential experimental design with
  applications in additive manufacturing.
\newblock {\em {arXiv} preprint {arXiv:2107.12809}}, July 2021.

\bibitem{goguelin2021}
S.~Goguelin, V.~Dhokia, and J.~M. Flynn.
\newblock Bayesian optimisation of part orientation in additive manufacturing.
\newblock {\em International Journal of Computer Integrated Manufacturing},
  34:1263--84, December 2021.

\bibitem{shapre2018}
{\em {Design of Mechanical Metamaterials via Constrained Bayesian
  Optimization}}, volume Volume 2A: 44th Design Automation Conference of {\em
  International Design Engineering Technical Conferences and Computers and
  Information in Engineering Conference}, 08 2018.
\newblock V02AT03A029.

\bibitem{xue2020}
T.~Xue, T.J. Wallin, Y.~Menguc, S.~Adriaenssens, and M.~Chiaramonte.
\newblock Machine learning generative models for automatic design of
  multi-material {3D} printed composite solids.
\newblock {\em Extreme Mechanics Letters}, 41:100992, November 2020.

\bibitem{hertlein2020}
N.~Hertlein, K.~Vemaganti, and S.~Anand.
\newblock Bayesian optimization of energy-absorbing lattice structures for
  additive manufacturing.
\newblock In {\em ASME International Mechanical Engineering Congress and
  Exposition}, volume 84539, 2020.

\bibitem{modal2020}
S.~Mondal, D.~Gwynn, A.~Ray, and A.~Basak.
\newblock Investigation of melt pool geometry control in additive manufacturing
  using hybrid modeling.
\newblock {\em Metals}, 10:683, May 2020.

\bibitem{xiong2019}
Y.~Xiong, P.~L. Duong, D.~Wang, S.~I. Park, Q.~Ge, N.~Raghavan, and D.~W.
  Rosen.
\newblock Data-driven design space exploration and exploitation for design for
  additive manufacturing.
\newblock {\em Journal of Mechanical Design}, 141, October 2019.

\bibitem{baturynska2018}
I.~Baturynska, O.~Semeniuta, and K.~Martinsen.
\newblock Optimization of process parameters for powder bed fusion additive
  manufacturing by combination of machine learning and finite element method: A
  conceptual framework.
\newblock In {\em Procedia Cirp}, volume~67, pages 227--32, January 2018.

\bibitem{govindjee}
S.~Reese and S.~Govindjee.
\newblock Theoretical and numerical aspects in the thermo-viscoelastic material
  behaviour of rubber-like polymers.
\newblock {\em Mechanics of Time-Dependent Materials}, 1:357--396, December
  1997.

\bibitem{pandey2003}
P.M. Pandey, N.V. Reddy, and S.G. Dhande.
\newblock Improvement of surface finish by staircase machining in fused
  deposition modeling.
\newblock {\em Journal of Materials Processing Technology}, 132:323--331,
  January 2003.

\bibitem{frazier2018tutorial}
Peter~I Frazier.
\newblock A tutorial on bayesian optimization.
\newblock {\em arXiv preprint arXiv:1807.02811}, 2018.

\bibitem{polak}
E.~Polak.
\newblock {\em Computational Methods in Optimization: A Unified Approach}.
\newblock Academic press, New York, 1971.

\bibitem{frazier2018}
P.~I. Frazier.
\newblock Bayesian optimization.
\newblock In {\em Recent advances in optimization and modeling of contemporary
  problems}, pages 255--278. Informs, October 2018.

\bibitem{shahriari2015taking}
Bobak Shahriari, Kevin Swersky, Ziyu Wang, Ryan~P Adams, and Nando De~Freitas.
\newblock Taking the human out of the loop: A review of bayesian optimization.
\newblock {\em Proceedings of the IEEE}, 104(1):148--175, 2015.

\bibitem{brochu2010}
E.~Brochu, V.~M. Cora, and N.~De Freitas.
\newblock A tutorial on bayesian optimization of expensive cost functions, with
  application to active user modeling and hierarchical reinforcement learning.
\newblock {\em {arXiv} preprint {arXiv:1012.2599}}, December 2010.

\bibitem{optimization}
J.~Nocedal and S.~Wright.
\newblock {\em Numerical Optimization}.
\newblock Springer-Verlag, New York, 2006.

\bibitem{scikit-learn}
F.~Pedregosa, G.~Varoquaux, A.~Gramfort, V.~Michel, B.~Thirion, O.~Grisel,
  M.~Blondel, P.~Prettenhofer, R.~Weiss, V.~Dubourg, J.~Vanderplas, A.~Passos,
  D.~Cournapeau, M.~Brucher, M.~Perrot, and E.~Duchesnay.
\newblock Scikit-learn: Machine learning in {P}ython.
\newblock {\em Journal of Machine Learning Research}, 12:2825--2830, 2011.

\bibitem{sklearn_api}
L.~Buitinck, G.~Louppe, M.~Blondel, F.~Pedregosa, A.~Mueller, O.~Grisel,
  V.~Niculae, P.~Prettenhofer, A.~Gramfort, J.~Grobler, R.~Layton,
  J.~VanderPlas, A.~Joly, B.~Holt, and G.~Varoquaux.
\newblock {API} design for machine learning software: experiences from the
  scikit-learn project.
\newblock In {\em ECML PKDD Workshop: Languages for Data Mining and Machine
  Learning}, pages 108--122, 2013.

\bibitem{choi2016influence}
Young-Hyu Choi, Cheol-Min Kim, Hwan-Seock Jeong, and Jeong-Ho Youn.
\newblock Influence of bed temperature on heat shrinkage shape error in fdm
  additive manufacturing of the abs-engineering plastic.
\newblock {\em World Journal of Engineering and Technology}, 4(3):186--192,
  2016.

\bibitem{paul2014effect}
Ratnadeep Paul, Sam Anand, and Frank Gerner.
\newblock {Effect of Thermal Deformation on Part Errors in Metal Powder Based
  Additive Manufacturing Processes}.
\newblock {\em Journal of Manufacturing Science and Engineering},
  136(3):031009, 03 2014.

\end{thebibliography}
\newpage
\clearpage
\dsp
\appendix
\normalsize
%
%
\setcounter{section}{1}
\setcounter{equation}{0}
\renewcommand{\thesection}{\Alph{section}}
\def\theequation{\Alph{section}.\arabic{equation}}
\addcontentsline{toc}{section}{Appendix A: Finite element matrices} 
\section*{Appendix: Finite element matrices}\label{sec:app}
\par\noindent
The spatial approximation of the dependent variables and their derivatives 
employs standard piecewise polynomials stemming from 4-node isoparametric 
element interpolation in two dimensions. 

The values of the material parameters are listed in 
Table~\ref{tbl:matcoeff-simp}. The material reference temperature is denoted by 
$\theta_{0,m}$. The parameters were estimated through least squares fitting of 
the experimental data for the temperature range of the simulations in this 
work. 
\begin{table}[H]
\centering
\begin{tabular}{|c|c|c|}
  \hline 
 Material parameter & Numerical values & Units\\
 \hline
$a$            & $-5.5$ & $-$ \\
$b$            & $1.0$  & $-$ \\
$\theta_{0,m}$ & $475$  & K\\
$E_{0,m}$      & $250$  &MPa\\
\hline
\end{tabular}
\caption{\em Material parameters for the simplified constitutive
function in Equation~\eqref{const-helm-energy-abs-simp}} 
\label{tbl:matcoeff-simp}
\end{table}
\par
The element-level matrices in~\eqref{eqn:disc-matrix-Newton} are 
\begin{equation}\label{eqn:Jacobian-matrix-Newton}
 \centering
 \begin{aligned}
   \left[\Mbm^e_{u}\right]\ =\ & \int_{\Omega_0^e}
\left[\Nbm_u^{e}\right]^T \rho_0 \left[\Nbm_u^{e}\right] \,dV\ ,\\
   \left[\Rbm^e_{u,n+1}\right]\ =\ &  \int_{\Omega_0^e} 
\left[\Bbm_u^{e}\right]^T \left<\sigmabold_{n+1}^e\right> \,dV,\\
   \left[\Fbm_{n+1}^e\right]\ =\ &\int_{\Omega_0^e}\left[\Nbm_u^{e}\right]^T \rho_0\left[\bbm_{n+1}\right] \,dV + 
                   \int_{\Gamma_{N,0}^u \cap
\partial\Omega_0^e}\left[\Nbm_u^{e}\right]^T
\left[\bar\pbm_{n+1}^e\right] \,dA\ ,\\
   \left[\Tbm^e_{n+1}\right]\ =\ & \int_{\Omega_0^e}
\left[\Nbm_t^{e}\right]^T c_{n+1}^e \left[\Nbm_t^{e}\right]\,dV\ ,\\
   \left[\Mbm^e_{t}\right]\ =\ & \int_{\Omega_0^e} k\left[\Bbm_t^{e}\right]^T\left[\Bbm_t^e\right] \,dV -  
     \int_{\Gamma_{N,0}^t\cap\partial\Omega_0^e}
\left[\Nbm_t^{e}\right]^T h \left[\Nbm_t^e\right]\,dA\ ,\\
   \left[\Rbm^e_{t,n+1}\right]\ =\ & \textcolor{black}{-}\int_{\Omega_0^e}\left[\Nbm_t^{e}\right]^T\theta_{n+1}^e 
    \left<\Mbm_{n+1}^e\right>^T \left<\dot{\epsbold}_{n+1}^e\right>\,dV\ ,\\
   \left[\Qbm^e_{n+1}\right]\ =\ & \int_{\Omega_0^e} \left[\Nbm_t^{e}\right]^T \rho_0 r_{n+1} \,dV -
      \int_{\Gamma_{N,0}^t\cap\partial\Omega_0^e}
\left[\Nbm_t^{e}\right]^T h \theta_{\infty}\,dA\ ,
  \end{aligned}
\end{equation}
and also 
\allowdisplaybreaks
\begin{align}\label{eqn:discrete-matrix-ele}
  \left[\Kbm^{e,(k)}_{u,n+1}\right]\ =\ & \frac{4}{\dt{n}^2} 
\int_{\Omega_0^e} \left[\Nbm_{u}^{e}\right]^T \rho_0 \left[\Nbm_{u}^{e}\right]\,dV 
+\int_{\Omega_0^e}\left[\Bbm_{u}^{e}\right]^T\left[\left(\pder{\sigmabold}{\epsbold}\right)_{n+1}^{e,(k)}\right]\left[\Bbm_{u}^{e}\right]\,dV\ ,\\
  \left[\Kbm^{e,(k)}_{t,n+1}\right]\ =\ & 
   \int_{\Omega_0^e} \left[\Bbm_{u}^{e}\right]^T
\left<\left(\pder{\sigmabold}{\theta}\right)_{n+1}^{e,(k)}\right>\left[\Nbm_{t}^{e}\right] \,dV\ ,\\
  \left[\Abm^{e,(k)}_{u,n+1}\right]\ =\ &
\int_{\Omega_0^e}\left[\Nbm_{t}^{e}\right]^T 
\left<\left(\pder{c}{\epsbold}\right)_{n+1}^{e,(k)}\right>
      \dot\theta_{n+1}^{e,(k)} \left[\Bbm_{u}^{e}\right] \,dV \\
     &-\textcolor{black}{\frac{2}{\dt{n}}}\int_{\Omega_0^e} \left[\Nbm_{t}^{e}\right]^T \theta_{n+1}^{e,(k)}\left<\Mbm_{n+1}^{e,(k)}\right>^T 
      \left[\Bbm_{u}^{e}\right] \,dV\\
     &-\int_{\Omega_0^e}  \left[\Nbm_{t}^{e}\right]^T\theta_{n+1}^{e,(k)}
\textcolor{black}{
\left<\dot{\epsbold}_{n+1}^{e,(k)}\right>^T 
\left[\left(\pder{\Mbm}
      {\epsbold}\right)_{n+1}^{e,(k)}\right]
}
\left[\Bbm_{u}^{e}\right] \,dV\ ,\\
  \left[\Abm^{e,(k)}_{t,n+1}\right]\ =\ 
    &\int_{\Omega_0^e} \left[\Nbm_{t}^{e}\right]^T \left[\left(\pder{c}{\theta}\right)_{n+1}^{e,(k)}\right]
      \dot\theta_{n+1}^{e,(k)} \left[\Nbm_{t}^{e}\right] \,dV \tag{\stepcounter{equation}\theequation} \\
    &+\int_{\Omega_0^e} \left[\Nbm_{t}^{e}\right]^T \frac{c_{n+1}^{e,(k)}}{\dt{n}} \left[\Nbm_{t}^{e}\right] \,dV
     +\int_{\Omega_0^e} k\left[\Bbm_{t}^{e}\right]^T \left[\Bbm_{t}^{e}\right] \,dV\\
    &-\int_{\Omega_0^e} \left[\Nbm_{t}^{e}\right]^T \left<\Mbm_{n+1}^{e,(k)}\right>^T\left<\dot{\epsbold}_{n+1}^{e,(k)}\right> \left[\Nbm_{t}^{e}\right]\,dV \\
    &-\int_{\Omega_0^e} \left[\Nbm_{t}^{e}\right]^T \theta_{n+1}^{e,(k)}\left<\left(\pder{\Mbm}
      {\theta}\right)_{n+1}^{e,(k)}\right>^T \left<\dot{\epsbold}_{n+1}^{e,(k)}\right>\left[\Nbm_{t}^{e}\right] \,dV \\
    &\textcolor{black}{-}\int_{\partial\Omega_0^e\cap\Gamma_{N,0}^t} \left[ \Nbm_{t}^{e}\right]^T  
     h \left[\Nbm_{t}^{e}\right] \,dA\ .
\end{align}
\clearpage
\par
To deduce the detailed forms of 
sensitivities, $\rbm_n^h$ is divided into discrete linear momentum and energy 
balances, denoted separately as~$\rbm_{1,n}^h$ and~$\rbm_{2,n}^h$, 
corresponding to two equations in~\eqref{eqn:discrete-balance-global}. 
Consequently, equation~(\ref{eqn:pde-general}) can be equivalently rewritten as 
\begin{equation}\label{eqn:pde-general-2}
 \centering
 \begin{aligned}
 \rbm_{1,n}^h \left( \ubm_{n}^h, \ubm_{n-1}^h, \vbm_{n-1}^h, \thetabold_{n}^h, \thetabold_{n-1}^h, \ubm^{his}_{n},
      \thetabold^{his}_{n},\ybm,t_{n}\right) \ =\ \zerobold,\\
 \rbm_{2,n}^h \left( \ubm_{n}^h, \ubm_{n-1}^h, \vbm_{n-1}^h, \thetabold_{n}^h, \thetabold_{n-1}^h, \ubm^{his}_{n},
      \thetabold^{his}_{n},\ybm,t_{n}\right) \ =\ \zerobold\ .
 \end{aligned}
\end{equation}
The element-level matrix-vector forms of~$\rbm_{1,n}^h$ and 
$\rbm_{2,n}^h$ are given in Equations 
~(\ref{eqn:discrete-balance-ele}) and (\ref{eqn:discrete-matrix-ele}).
Thus, excluding any interelement contributions, the element-level matrices of
$\pder{\rbm_{n}^h}{\ubm_{n-1}^h}, \ \pder{\rbm_{n}^h}{\vbm_{n-1}^h}, \ \pder{\rbm_{n}^h}{\thetabold_{n-1}^h}$ are 
\begin{equation}\label{eqn:drdun}
 \centering
 \begin{aligned}
 &\left[\pder{\rbm_{1,n}^{e,h}}{\ubm_{n-1}^h}\right]\ =\ -\frac{4}{\dt{n}^2}\left[\Mbm_{u}^e\right], \qquad
 &&\left[\pder{\rbm_{2,n}^{e,h}}{\ubm_{n-1}^h}\right]\ =\  \int_{\Omega_0^e} \left[\Nbm_{t}^{e}\right]^T \left[\theta_{n+1}^e\right] \left[\Mbm_{n+1}^e\right] \frac{1}{\dt{n}} 
                         \left[\Bbm_{u2}^{e}\right] ~dV\ ,\\
 &\left[\pder{\rbm_{1,n}^{e,h}}{\vbm_{n-1}^h}\right]\ =\  -\frac{4}{\dt{n}}\left[\Mbm_{u}^e\right], \qquad
 &&\left[\pder{\rbm_{2,n}^{e,h}}{\vbm_{n-1}^h}\right]\ =\  \int_{\Omega_0^e} \left[\Nbm_{t}^{e}\right]^T \theta_{n+1}^e \left[\Mbm_{n+1}^e\right]
                         \left[\Bbm_{u2}^{e}\right] ~dV\ ,\\
 &\left[\pder{\rbm_{1,n}^{e,h}}{\thetabold_{n-1}^h}\right]\ =\ \quad  \left[\zerobold\right], \qquad
 &&\left[\pder{\rbm_{2,n}^{e,h}}{\thetabold_{n-1}^h}\right]\ =\   
     -\int_{\Omega_0^e} \left[\Nbm_{t}^{e}\right]^T
\frac{c_{n+1}^k}{\dt{n}} \left[\Nbm_{t}^{e}\right] ~dV\ .
 \end{aligned}
\end{equation}
Assemble Equations~(\ref{eqn:drdun}) over the entire domain. The global 
matrix forms of the sensitivities terms 
\begin{equation}\label{pde-diff-disc-1}
 \centering
   \left[\pder{\rbm_{i,n}^h}{\ubm_{n-1}^h}\right]\ =\ \assembly \left[\pder{\rbm_{i,n}^{e,h}}{\ubm_{n-1}^h}\right],
    \left[\pder{\rbm_{i,n}^h}{\vbm_{n-1}^h}\right]\ =\ \assembly\left[\pder{\rbm_{i,n}^{e,h}}{\vbm_{n-1}^h}\right] ,
   \left[\pder{\rbm_{n}^h}{\thetabold_{n-1}^h}\right]\ =\
\assembly\left[\pder{\rbm_{2,n}^{e,h}}{\thetabold_{n-1}^h}\right]\ ,
\end{equation}
$i=1,2$ are obtained.
At the element level, 
\begin{equation}\label{eqn:sens-6}
 \centering
\pder{\rbm_{i,n}^{e,h}}{\ubm^{his}_{n}},  \pder{\rbm_{i,n}^{e,h}}{\thetabold^{his}_{n}}, \ i=1,2, 
\end{equation}
can be calculated similarly in matrix form by differentiating the fully 
discretized form in Equations~(\ref{eqn:discrete-balance-ele}) 
and~(\ref{eqn:discrete-matrix-ele}): 
\begin{equation}\label{eqn:drdhis}
 \centering
\begin{aligned}
\left[\pder{\rbm_{1,n}^{e,h}}{\ubm^{his}_{n}}\right] \ =\ &
-\int_{\Omega_0^e}\left[\Bbm_{u}^{e}\right]^T\left[\left(\pder{\sigmabold}{\epsbold}\right)_{n+1}^{e}\right]\left[\Bbm_{u2}^{e}\right]~dV\ ,\\
\left[\pder{\rbm_{1,n}^{e,h}}{\thetabold^{his}_{n}}\right] \
=\ &  \int_{\Omega_0^e} \left[\Bbm_{u}^{e}\right]^T 
           \left[\left(\pder{\sigmabold}{\thetabold^{his}}\right)_{n}^{e}\right]\left[\Nbm_{t}^{e}\right] ~dV\ ,\\
\left[\pder{\rbm_{2,n}^{e,h}}{\ubm^{his}_{n}}\right] \ =\ &  \int_{\Omega_0^e} \left[\Nbm_{t}^{e}\right]^T \left[\left(\pder{c}{\ubm^{his}}\right)_{n+1}^e\right]
      \dot\theta_{n}^e \left[\Bbm_{u2}^{e}\right] ~dV -\\
    &\int_{\Omega_0^e} \left[\Nbm_{t}^{e}\right]^T\theta_{n}^e\left[\left(\pder{\Mbm}
      {\epsbold}\right)_{n}^e\right]
\left[\dot{\epsbold}_{n}^e\right]\left[ \Bbm_{u2}^{e}\right] ~dV\ , \\ 
\left[\pder{\rbm_{2,n}^{e,h}}{\thetabold^{his}_{n}}\right]\ =\ 
&\int_{\Omega_0^e} \left[\Nbm_{t}^{e}\right]^T \left[\left(\pder{c}{\theta^{his}}\right)_{n}^{e}\right]
      \dot\theta_{n}^{e,(k)} \left[\Nbm_{t}^{e}\right] \,dV  \\
    &-\int_{\Omega_0^e} \left[\Nbm_{t}^{e}\right]^T \theta_{n}^{e}\left<\left(\pder{\Mbm}
{\theta^{his}}\right)_{n}^{e}\right>^T \left<\dot{\epsbold}_{n}^{e}\right>\left[\Nbm_{t}^{e}\right]~dV. \\
\end{aligned}
\end{equation}
\\par
This primal sensitivity framework can be applied potentially to any time-step 
independent and differentiable process parameters. As an example, 
if the framework is applied to the convection coefficient~$h$, 
which can be considered a constant, the corresponding element-level 
sensitivity of~\eqref{eqn:pde-diff} is 
\begin{equation}\label{eqn:dfdch}
 \centering
\begin{aligned}
\left[\pder{\rbm_{1,n}^{e,h}}{h}\right] \ =\ & \quad \left[\zerobold\right], \\
\left[\pder{\rbm_{2,n}^{e,h}}{h}\right]\ =\ &
\int_{\Gamma_{N,0}^t\cap\partial\Omega_0^e}  \left[\Nbm_t^{e}\right]^T
\left(\left[\thetabold_{\infty}\right]-\left[\thetabold_{n}^e\right]
\right)~dA\ .
\end{aligned}
\end{equation}
Assembling Equation~(\ref{eqn:dfdch}) over the entire domain 
and the matrix form of $\pder{\rbm_{n}^h}{\ybm}$ is obtained. 
%
%
%
\clearpage
\end{document}